Гасников Александр Владимирович

# Эффективные численные методы поиска равновесий в больших транспортных сетях

Специальность 05.13.18 – Математическое моделирование, численные методы и комплексы программ

Автореферат диссертации на соискание ученой степени доктора физико-математических наук

МОСКВА – 2016

Работа выполнена на кафедре математических основ управления факультета управления и прикладной математики (ФУПМ) Московского физико-технического института (государственного университета)

**Научный консультант:**

доктор физ.-мат. наук, профессор Шананин Александр Алексеевич

**Официальные оппоненты:**

доктор техн. наук, старший научный сотрудник Поляк Борис Теодорович
Институт проблем управления им. В.А. Трапезникова,
Лаборатория № 7 им. Я.З. Цыпкина, главный научный сотрудник

доктор физ.-мат. наук, доцент Посыпкин Михаил Анатольевич
Федеральный Исследовательский Центр "Информатика и Управление" РАН, Вычислительный центр РАН, Отдел прикладных проблем оптимизации, и.о. заведующего отделом

доктор физ.-мат. наук, профессор Васин Александр Алексеевич
кафедра исследования операций факультета Вычислительной математики и кибернетики МГУ им. М.В. Ломоносова, и.о. заведующего кафедрой

**Ведущая организация:**

Центральный экономико-математический институт РАН (г. Москва)

Защита состоится 1 декабря 2016 г. в 10 ч. 45 мин. на заседании диссертационного совета Д 212.156.05 на базе Московского физико-технического института (государственного университета) по адресу: 141700, Московская обл., г. Долгопрудный, Институтский пер., д. 9, ауд. 903 КПМ.

С диссертацией можно ознакомиться в библиотеке МФТИ (ГУ) и на сайте МФТИ https://mipt.ru/

Автореферат разослан 14 сентября 2016 года.

Ученый секретарь дис. совета                                          Федько Ольга Сергеевна

# Общая характеристика работы

## Актуальность темы и степень ее разработанности

Настоящая диссертация посвящена разработке новых подходов к построению многостадийных моделей транспортных потоков и эффективных численных методов поиска равновесий в таких моделях. Начиная с 50-х годов XX века вопросам поиска равновесий в транспортных сетях стали уделять большое внимание в связи с ростом городов и необходимостью соответствующего транспортного планирования. В 1955 г. появилась первая модель равновесного распределения потоков по путям: BMW-модель, также называемая моделью Бэкмана[1]. В этой модели при заданных корреспонденциях (потоках из источников в стоки) решалась задача поиска равновесного распределения этих корреспонденций по путям, исходя из принципа Вардропа, т.е. исходя из предположения о том, что каждый пользователь транспортной сети рационален и выбирает кратчайший маршрут следования. Таким образом, поиск равновесия в такой модели сводился к поиску равновесия Нэша в популяционной игре[2] (популяций столько, сколько корреспонденций). Поскольку в модели предполагалось, что время прохождения ребра есть функция от величины потока только по этому ребру, то получившаяся игра была игрой загрузки, следовательно, потенциальной. Последнее означает, что поиск равновесия сводится к решению задачи оптимизации. Получившуюся задачу выпуклой оптимизации решали с помощью метода условного градиента[3]. Описанная модель и численный метод и по настоящее время используются в подавляющем большинстве продуктов транспортного моделирования для описания блока равновесного распределения потоков по путям[4]. Однако в работе[5] Нестерова–де Пальмы было указано на ряд существенных недостатков модели Бэкмана, и предложена альтернативная модель, которую авторы назвали моделью Стабильной Динамики.

---

[1] Beckmann M., McGuire C.B., Winsten C.B. Studies in the economics of transportation. RM-1488. Santa Monica: RAND Corporation, 1955.

[2] Sandholm W. Population games and Evolutionary dynamics. Economic Learning and Social Evolution. MIT Press; Cambridge, 2010.

[3] Frank M., Wolfe P. An algorithm for quadratic programming // Naval research logistics quarterly. – 1956. V. 3. № 1-2. P. 95–110.

[4] Стенбринк П.А. Оптимизация транспортных сетей. М.: Транспорт, 1981; Sheffi Y. Urban transportation networks: Equilibrium analysis with mathematical programming methods. N.J.: Prentice–Hall Inc., Englewood Cliffs, 1985; Patriksson M. The traffic assignment problem. Models and methods. Utrecht, Netherlands: VSP, 1994. Ortúzar J.D., Willumsen L.G. Modelling transport. JohnWilley & Sons, 2011.

[5] Nesterov Y., de Palma A. Stationary Dynamic Solutions in Congested Transportation Networks: Summary and Perspectives // Networks Spatial Econ. – 2003. – № 3(3). – P. 371–395.

Отмеченные выше модели равновесного распределения потоков по путям могут быть использованы при решении различных задач долгосрочного планирования. Например, такой задачи. Имея заданный бюджет, нужно решить, на каких участках графа транспортной сети стоит увеличить полосность дороги / построить новые дороги. Заданы несколько сценариев, нужно отобрать лучший. Задачу можно решить, найдя равновесные распределения потоков, отвечающие каждому из сценариев, и сравнивая найденные решения, например, по критерию суммарного времени, потерянного в пути всеми пользователями сети в данном равновесии. При значительных изменениях графа транспортной сети необходимо в приведенную выше цепочку рассуждений включать дополнительный контур, связанный с тем, что изменения приведут не только к перераспределению потоков на путях, но и поменяют корреспонденции. Таким образом, корреспонденции также необходимо моделировать. В 60-е годы XX века появилось сразу несколько различных моделей для расчета матрицы корреспонденций, исходя из информации о численностях районов и числе рабочих мест в них. Наибольшую популярность приобрела энтропийная модель расчета матрицы корреспонденций (А.Дж. Вильсон). В этой модели поиск матрицы корреспонденций сводился к решению задачи энтропийно-линейного программирования.

К сожалению, при этом в энтропийную модель явным образом входит информация о матрице затрат на кратчайших путях по всевозможным парам районов. Возникает "порочный круг": чтобы посчитать эту матрицу затрат, нужно сначала найти равновесное распределение потоков по путям, а чтобы найти последнее, необходимо знать матрицу корреспонденций, которая рассчитывается по матрице затрат. На практике отмеченную проблему решали методом простых итераций. Как-то "разумно" задавали начальную матрицу корреспонденций, по ней считали распределение потоков по путям, по этому распределению считали матрицу затрат, на основе которой пересчитывали матрицу корреспонденций, и процесс повторялся. Повторялся он до тех пор, пока не выполнялся критерий останова. К сожалению, до сих пор для описанной процедуры не известно никаких гарантий ее сходимости и, тем более, оценок скорости сходимости.

Описанный выше подход можно назвать двухстадийной моделью транспортных потоков, потому что модель состоит из прогонки двух разных блоков. В действительности, в реальных приложениях число блоков 3-4. В част-ности, как правило, всегда включают блок расщепления потоков по типу передвижения (например, личный и общественный транспорт) – этот блок описывается моделью аналогичной модели Бэкмана. В математическом плане это уточнение не существенно. Все основные имеющиеся тут сложности хорошо демонстрирует уже двухстадийная модель. Отметим также, что часто в приложениях вместо модели Бэкмана используется её "энтропийно-регуляризованный" вариант, который отражает ограниченную рациональность пользователей транспортной сети. Равновесие в такой модели часто называют "стохастическим" равновесием.

Резюмируем написанное выше. До настоящего момента не существовало строгого научного обоснования используемого повсеместно на практике (и зашитого во все современные пакеты транспортного моделирования) способа формирования многостадийных моделей транспортных потоков. Не существовало также никаких гарантий сходимости численных методов, используемых для поиска равновесий в многостадийных моделях. В используемых сейчас повсеместно многостадийных моделях в качестве основных блоков фигурируют блоки с моделями типа Бэкмана, а не более современные блоки Стабильной Динамики. Таким образом, актуальной является задача обоснования современной многостадийной модели и разработка эффективных численных методов поиска (стохастического) равновесия в такой модели.

**Цели и задачи**

Многие законы природы могут быть записаны в форме "вариационных принципов". В моделировании транспортных потоков это также имеет место. Однако, на текущий момент с помощью вариационных принципов описываются только отдельные блоки многостадийной транспортной модели, и эволюционный вывод вариационных принципов имеется только для блоков с моделью Бэкмана в основе. Одной из целей данной диссертационной работы является эволюционный вывод всех блоков многостадийной транспортной модели (прежде всего, речь идет о блоке расчета матрицы корреспонденций), и получение (с помощью эволюционного вывода) вариационного принципа для описания равновесия в многостадийной модели.

Целью также является разработка "алгебры" над блоками-моделями (каждый блок описывается своим вариационным принципом), которая позволит, как в конструкторе, собирать (формируя общий вариационный принцип) сколько угодно сложные модели из небольшого числа базисных элементов конструктора (блоков).

Описанный выше формализм приводит в итоге к решению задач выпуклой оптимизации в пространствах огромных размеров и имеющих довольно специальную иерархическую (многоуровневую) структуру функционала задачи. Чтобы подчеркнуть нетривиальность таких задач, отметим, что переменные, по которым необходимо оптимизировать, – это, в частности, компоненты вектора распределения потоков по путям. Для графа в виде двумерной квадратной решетки (Манхэттенская сеть) с числом вершин порядка нескольких десятков тысяч (это число соответствует транспортному графу Москвы) такой вектор с большим запасом нельзя загрузить в память любого современного суперкомпьютера, не говоря уже о том, чтобы как-то работать с такими векторами.

Важной целью диссертации является разработка (с теоретическими гарантиями) эффективных численных методов, способных за несколько часов на персональном компьютере с хорошей точностью (и с высокой вероятностью) найти равновесие в многостадийной модели транспортных потоков крупного мегаполиса.

В частности, целью является разработка "алгебры" над численными методами, используемыми для расчета отдельных блоков многостадийной модели, которая позволит, как в конструкторе, собирать итоговую эффективную численную процедуру (для поиска равновесия в многостадийной модели) с помощью правильного чередования / комбинации работы алгоритмов, используемых для отдельных блоков.

Более общей целью является выделение небольшого набора "оптимальных" базисных численных методов выпуклой оптимизации и операций над ними, чтобы с помощью всевозможных сочетаний можно было получать "оптимальные" методы для структурно сложных задач выпуклой оптимизации. То есть цель – научиться раскладывать (декомпозировать) сложно составленную задачу выпуклой оптимизации на простые блоки, чтобы численный метод для общей задачи можно было бы "собрать" из простых базисных бло-ков (методов). Желательно также, чтобы разработанный формализм позволял автоматизировать эту процедуру.

**Научная новизна, методология и методы исследования**

В основе предложенного в диссертации эволюционного формализма обоснования многостадийной транспортной модели лежит часто используемая в популяционной теории игр марковская logit-динамика, отражающая ограниченную рациональность агентов (водителей). Новым является понимание этой динамики как модели стохастической химической кинетики с унарными реакциями и рассмотрение сразу нескольких разных типов таких унарных реакций, происходящих с разной (по порядку величины) интенсивностью, и отвечающих разномасштабным процессам, протекающим в транспортной системе. Например, для двухстадийной модели динамика, отвечающая формированию корреспонденций, идёт, по терминологии А.Н. Тихонова, в медленном времени (годы), а отвечающая распределению потоков по путям – в быстром времени (дни). Тогда с некоторыми оговорками функционал в вариационном принципе с точностью до множителя и аддитивной константы можно, с одной стороны, понимать как функционал Санова (действие), отвечающий за концентрацию стационарной (инвариантной) меры введенной марковской динамики, а, с другой стороны, как функционал Ляпунова–Больцмана для кинетической динамики, полученной при (каноническом) скейлинге (по числу агентов) введённой марковской динамики.

При разработке эффективных численных процедур в качества базиса были выбраны два метода: метод зеркального спуска – МЗС (А.С. Немировский, 1977) и быстрый градиентный метод – БГМ (Ю.Е. Нестеров, 1982). Практически все используемые в диссертации алгоритмы являются некоторыми вариантами (производными) этих двух методов.

В диссертации часто использовались следующие операции над алгоритмами: процедура рестартов, регуляризация функционала задачи, mini-batc'инг. По отдельности эти процедуры были давно и хорошо известны. Однако в работе были предложены различные сочетания отмеченных операций, позволившие получить часть результатов.

Важное место в диссертации занимает "игра" между стоимостью итерации численного метода и необходимым числом итераций. Здесь было рассмотрено два направления.

1. В большинстве приложений "стоимость" (время) получения от оракула (роль которого, как правило, играют нами же написанные подпрограммы вычисления градиента) градиента функционала заметно превышает время, затрачиваемое на то, чтобы сделать шаг итерации, исходя из выданного оракулом вектора. Желание сбалансировать это рассогласование (усложнить итерации, сохранив при этом старый порядок их сложности, и выиграть за счет этого в сокращении числа итераций) привело к возникновению композитной оптимизации (Ю.Е. Нестеров, 2007), в которой (аддитивная) часть функционала задачи переносится без лианеризации (запроса градиента) в итерации. Здесь остается еще много степеней свободы, позволяющих играть на том, насколько "дорогим" окажется оракул и соответствующая этому оракулу "процедура проектирования", и на том, сколько (внешних) итераций потребуется методу для достижения заданной точности. В частности, если обращение к оракулу за градиентом и последующее проектирование требуют, в свою очередь, решения вспомогательных оптимизационных задач, то можно "сыграть" на том, насколько точно надо решать эти вспомогательные задачи, пытаясь найти компромисс между "стоимостью" итерации и числом итераций. Также можно сыграть и на том, как выделять эти вспомогательные задачи – иными словами, что понимать под оракулом и что – под итерацией метода. Общая идея подхода "разделяй и властвуй" к численным методам выпуклой оптимизации может принимать довольно неожиданные и при этом весьма эффективные формы (как например, в методах внутренней точки Нестерова–Немировского, 1989). В диссертации рассмотрены разные оригинальные варианты описанной игры в связи с транспортно-сетевыми приложениями. Заметим, что для обоснования упомянутых конструкций (подходов) в диссертации существенным образом используется (и развивается) концепция неточного оракула, выдающего градиент и значение функции (Деволдер–Глинёр–Нестеров, 2013). Эта же концепция существенным образом используется (и развивается) в диссертации при разработке линейки универсальных методов (Ю.Е. Нестеров, 2013), самонастраивающихся на гладкость задачи. Общая идея таких методов – искусственно вводить неточность, чтобы правильно подбирать класс гладкости задачи.

2. Достаточно часто имеет смысл (с точки зрения минимизации общего времени работы метода) вместо градиента в итерационную процедуру подставлять некоторый (легко вычислимый) аналог градиента. Как правило, это несмещенная (или мало смещенная) оценка градиента (стохастический градиент). Число итераций при этом возрастает, но зато итерации становятся более "дешевыми". В диссертации описываются общие "рецепты" получения из детерминированных (полноградиентных) методов их стохастических (рандомизированных) вариантов, позволяющие достаточно просто изучать "наследуемые" при этом свойства исходных методов. К таким рандомизированным методам можно отнести, например, спуски по направлению, покомпонентные методы и методы нулевого порядка. Нетривиально уже то, что это оказалось возможным. Например, покомпонентный БГМ был предложен (Ю.Е. Нестеров, 2012) практически независимо от полноградиентного варианта БГМ. В диссертации продемонстрировано, как можно достаточно просто получить покомпонентный БГМ с наследованием всех основных свойств из БГМ в специальной форме Allen-Zhu–Orecchia. Полезно заметить, что в такой форме БГМ можно проинтерпретировать, как выпуклую комбинацию МЗС и обычного метода проекции градиента. Важное место в диссертации занимает изучение описанного в этом пункте формализма вместе с концепцией эффективного пересчета используемого варианта градиента. Не расчета, а именно пересчета, т.е. расчета с учетом результатов вычислений, проделанных на предыдущих итерациях. Поясним последнее примером. При безусловной минимизации квадратичной формы с помощью (неускоренного) покомпонентного спуска новая точка итерационного процесса отличается от старой только в одной компоненте, поэтому можно так организовать вычисления, чтобы итерация в среднем занимала $O(s)$ арифметических операций ($s$ – среднее число ненулевых элементов в столбце матрицы квадратичной формы), независимо от размеров матрицы квадратичной формы. Развивая идеи Ю.Е. Нестерова, в рамках описанного выше формализма в диссертации изучены способы наилучшего учета разреженности задачи.

Ранее уже отмечалась проблема огромной размерности пространства потоков по путям, в котором ставятся задачи выпуклой оптимизации, возникающие при поиске равновесий в транспортных сетях. С помощью метода условного градиента эта проблема решается благодаря эффективно вычислимому с помощью алгоритма Дейкстры поиска кратчайших путей в графе линейному минимизационному оракулу. В диссертации выбран иной форма-лизм, более удобный для перенесения на многостадийные модели, связанный с переходом к двойственной задаче, и ее решением прямо-двойственным методом, позволяющим практически "бесплатно" восстанавливать решение прямой задачи. Отметим, что оба базисных метода (МЗС и БГМ) – прямо-двойственные. Этот формализм распространяется в диссертации на многостадийные модели. Именно с помощью переходов к двойственным задачам в части блоков удалось свести поиск равновесия в многостадийной модели к эффективно решаемой (с помощью описанных выше конструкций) задаче выпуклой оптимизации. Однако все это потребовало серьезного погружения в прямо-двойственность численных методов выпуклой оптимизации, особенно для задач выпуклой оптимизации на неограниченных областях. В частности, потребовалось изучение сочетания в одной задаче на одном методе основных (базисных) прямо-двойственных конструкций (способов восстановления решения сопряженной задачи), используемых ранее только по отдельности (Ю.Е. Нестеров, 2009; А.С. Немировский, 2010).

**Теоретическая и практическая значимость работы**

Настоящая диссертация мотивирована транспортными приложениями. Выбор задач был обусловлен общением со специалистами из НИиПИ Генплана г. Москвы, Департамента транспорта г. Москвы, компании А+С (основным дистрибьютором в России пакетов транспортного моделирования линейки PTV), ЦИТИ г. Москвы, Института экономики транспорта и транспортной политики НИУ ВШЭ.

За последние 10 лет резко возрос объем и качество доступных для моделирования транспортных данных (GPS-треки, данные сотовых операторов, данные видеокамер, всевозможные опросы населения, данные геоинформационных систем). При этом задачи поиска равновесий в транспортных сетях "вышли на передний план", поскольку появилась возможность ставить такие задачи на больших (детализированных) транспортных сетях. В результате сложность таких задач резко возросла. С другой стороны, увеличилась и потребность в многократном решении таких задач с целью просмотра различных сценариев развития транспортной инфраструктуры. Как следствие, появилась необходимость в разработке нового (адаптированного под эти реалии) аппарата моделирования и соответствующих вычислительных процедур.

Именно этому – поиску равновесий в больших (реальных) транспортных сетях, прежде всего, и посвящена диссертационная работа.

Большой акцент в диссертации сделан на теоретическое исследование оптимальных вычислительных процедур. Основные подходы здесь были заложены в Советском Союзе в работах Б.Т. Поляка, А.С. Немировского, Ю.Е. Нестерова и др.[6] В диссертации удалось "овыпуклить" многие классические результаты и посмотреть на многообразие этих результатов с единых позиций: с помощью базисного набора методов и введенных операций над ними удалось получить более простым способом как известные результаты, так и новые. В частности, исследовать вопросы: о равномерной ограниченности последовательностей, генерируемых численными методами (в том числе рандомизированными); о практической реализации всех рассматриваемых методов (с теоретическим обоснованием) в условиях отсутствия априорной информации о свойствах функционала задачи и свойствах решения (в частности, ограниченности нормы решения известным числом); о практически эффективных критериях останова рассматриваемых методов (с установленными теоретическими гарантиями скорости сходимости, что особенно нетривиально для прямо-двойственных методов на неограниченных областях).

Отметим также, что подавляющая часть изложения в диссертации ведется на современном уровне[7] – с точными константами в оценках числа итераций и с оценками вероятностей больших уклонений для стохастических (рандомизированных) методов.

Алгоритмы поиска равновесий в транспортных сетях, разработанные в диссертации, вошли в комплекс программ, созданных коллективом при участии автора диссертации. Комплекс был успешно принят в рамках отчета по гранту федеральной целевой программы «Исследования и разработки по приоритетным направлениям развития научно-технологического комплекса России на 2014 – 2020 годы», Соглашение № 14.604.21.0052 от 30.06.2014 г. с Минобрнаукой (уникальный идентификатор проекта RFMEFI60414X0052).

---

[6] Немировский А.С., Юдин Д.Б. Сложность задач и эффективность методов оптимизации. М.: Наука, 1979. – 384 с.; Поляк Б.Т. Введение в оптимизацию. М.: Наука, 1983. – 384 с.

[7] Нестеров Ю.Е. Введение в выпуклую оптимизацию. М.: МЦНМО, 2010. – 280 стр.; Nemirovski A. Lectures on modern convex optimization analysis, algorithms, and engineering applications. Philadelphia: SIAM, 2013. – 590 p.; Bubeck S. Convex optimization: algorithms and complexity // Foundations and Trends in Machine Learning. – 2015. – V. 8. no. 3-4. – P. 231–357.

Также на основе разработок главы 3 диссертации, в том числе автором диссертации, был создан комплекс программ, проданный компании Huawei Russia (Договор Huawei с МФТИ №: YB2014120038 от 23 декабря 2014 года).

**Положения, выносимые на защиту (основные положения – пп. 1 – 4)**

1. Предложен новый эволюционный вывод энтропийной модели расчета матрицы корреспонденций с помощью марковской logit-динамики и транспортных потенциалов Канторовича–Гавурина. Таким образом, показано, что на модель расчета матрицы корреспонденций можно смотреть как на (энтропийно-регуляризованную) разновидность модели Бэкмана.

2. Получен оригинальный вывод модели Стабильной Динамики из модели Бэкмана с помощью вырождения функций затрат на прохождения ребер методом внутренних штрафов.

3. Используя то, что все блоки многостадийной модели транспортных потоков есть вариации модели Бэкмана, являющейся, в свою очередь, популяционной игрой загрузки, предложен общий способ формирования вариационных принципов для поиска равновесий в многостадийных транспортных моделях. Полученные результаты распространены на общие иерархические популяционные игры загрузки. Исследован класс сетевых рынков (частным случаем которых является модель грузоперевозок РЖД) поиск конкурентных равновесий Вальраса в которых может быть осуществлен в рамках описанного выше (вариационного) формализма с заменой задачи выпуклой оптимизации на задачу поиска седловой точки с правильной выпукло-вогнутой структурой.

4. Впервые предложена вариация универсального быстрого градиентного метода Ю.Е. Нестерова (БГМ), самонастраивающегося на гладкость задачи для сильно выпуклых задач композитной оптимизации. Предложенную вариацию удалось распространить и на задачи стохастической оптимизации. Ранее считалось, что такая вариация либо невозможна, либо будет весьма сложной – полученные результаты опровергли эти опасения. Получены оценки скорости сходимости, показывающие, что предложенный метод является равномерно оптимальным (по числу итераций, с точностью до числового множителя) для общего класса задач выпуклой оп-тимизации. Описываемая линейка универсальных методов активно использовалась в диссертации при разработке комплекса программ для расчета различных блоков многостадийной транспортной модели.

5. Для поиска равновесия в модели Стабильной Динамики предложена специальная рандомизированная версия прямо-двойственного метода зеркального спуска (МЗС). Впервые получены оценки скорости сходимости в терминах вероятностей больших уклонений общего стохастического варианта МЗС для случая, когда оптимизация проводится на неограниченном множестве (как в модели Стабильной Динамики). Также впервые получены оценки вероятностной локализации итерационной последовательности. Используемая при этом оригинальная техника априорной формулировки гипотезы о характере хвостов распределений с апостериорной проверкой представляется полезной и для ряда других приложений.

6. Предложены два различных прямо-двойственных подхода для решения задачи энтропийно-линейного программирования, в частности, для задачи расчета матрицы корреспонденций. В основе обоих подходов лежит идея решения двойственной задачи с помощью БГМ. Оба подхода распространены на общие задачи минимизации сильно выпуклых функционалов простой (например, сепарабельной) структуры при аффинных ограничениях.

7. Исследованы общие способы приближенного восстановления решения сопряженной (двойственной) задачи по последовательности, генерируемой методом, решающим исходную задачу. В частности, было продемонстрировано сочетание различных способов восстановления на одной задаче – поиск равновесного распределения потоков по ребрам в Смешанной Модели, когда часть ребер – из модели Бэкмана, оставшаяся часть – из модели Стабильной Динамики.

8. Исследована роль неточностей неслучайной природы при расчете градиента и значения оптимизируемой функции, а также неточностей, возникающих при проектировании, на итоговые оценки скорости сходимости различных методов. Полученные здесь результаты используются для обоснования концепции суперпозиции (правильного чередования – чередования с правильными частотами) численных методов.

9. Разработана концепция суперпозиции численных методов выпуклой оптимизации. Рассмотрены конкретные примеры приложений (min max-задачи, min min-задачи), в частности, задачи поиска равновесия в многостадийных моделях транспортных потоков.

10. Исследованы покомпонентные методы, спуски по направлению и методы нулевого порядка, с точки зрения наследования "хороших" свойств своих полноградиентных аналогов. В частности, была исследована прямо-двойственность покомпонентных методов. Новые результаты удалось получить для спусков по направлению и их дискретных аналогов (методов нулевого порядка) в части сочетания структуры множества, на котором происходит оптимизация, и способа выбора случайного направления спуска. В зависимости от контекста, оказалось, что оптимально выбирать случайное направление либо среди координатных осей, либо равномерно на евклидовой сфере. Исследована возможность использования прямо-двойственных блочно-покомпонентных методов для поиска равновесия в модели Стабильной Динамики.

11. На примере задачи ранжирования web-страниц (Google problem) исследованы основные конструкции huge-scale-оптимизации. В частности, предложено учитывать разреженность задачи и использовать рандомизированные методы. Предложенные методы позволили не только эффективно находить вектор PageRank, но и решать целый ряд других задач huge-scale-оптимизации. Например, предложен эффективный численный метод решения равномерно разреженных (по строкам и столбцам) систем линейных уравнений огромных размеров. Доказана высокая эффективность метода, если решение состоит из разных по масштабу компонент (1-норма вектора решения и его 2-норма достаточно близки).

12. Алгоритм Григориадиса–Хачияна поиска равновесий в антагонистических матричных играх удалось проинтерпретировать как специальным образом рандомизированный МЗС, что позволило получить для него оценки вероятностей больших отклонений и распространить его на задачи с равномерно разреженной матрицей, показав, что при этом разреженность учитывается оптимальным образом. Также получена естественная онлайн-интерпретация этого алгоритма.

**Степень достоверности, публикации и апробация результатов**

Степень достоверности и обоснованности полученных результатов достаточно высоки. Они обеспечиваются строгостью и корректным использованием математических доказательств, подтверждением результатов работы экспертными оценками специалистов. Результаты, включенные в данную работу, представлены более чем в 40 работах, в числе которых одна монография, одно учебное пособие и 26 статей в журналах из перечня ВАК. В работах, написанных в соавторстве, соискателю принадлежат результаты, указанные в положениях, выносимых на защиту. По результатам диссертации было прочитано несколько оригинальных курсов лекций студентам МФТИ, НМУ, БФУ им. И. Канта. Большая часть этих лекций (более 50) доступна в сети Интернет: http://www.mathnet.ru/ и https://www.youtube.com/user/PreMoLab.

Результаты, включенные в диссертацию, неоднократно докладывались в 2010–2016 гг. на научных семинарах в ведущих университетах и институтах России: семинаре отдела Математического моделирования экономических систем ВЦ ФИЦ ИУ РАН (рук. чл.-корр. РАН И.Г. Поспелов), семинаре лаборатории Адаптивных и робастных систем им. Я.З. Цыпкина ИПУ РАН (рук. проф. Б.Т. Поляк), семинаре лаборатории Структурных методов анализа данных в предсказательном моделировании МФТИ–ИППИ (рук. проф. В.Г. Спокойный), семинаре лаборатории Больших случайных систем мехмата МГУ (рук. проф. В.А. Малышев), семинаре кафедры Исследования операций ВМиК МГУ (рук. проф. А.А. Васин), семинаре ЦЭМИ РАН (рук. проф. Л.А. Бекларян), семинаре ИАП РАН (рук. чл.-корр. РАН А.С. Холодов), семинаре ПОМИ РАН (рук. проф. А.М. Вершик), семинарах федерального профессора А.М. Райгородского в МФТИ и МГУ, семинарах Яндекса (Москва, Новосибирск), семинаре кафедры Высшей математики МФТИ (рук. проф. Е.С. Половинкин), семинаре кафедры Информатики МФТИ (рук. чл.-корр. РАН И.Б. Петров), семинаре НИиПИ Генплана г. Москвы, семинаре по Транспортному моделированию ИПМ РАН (рук. доц. В.П. Осипов, акад. Б.Н. Четверушкин), семинаре МАДИ–МТУСИ–МИАН (рук. проф. А.П. Буслаев, проф. М.В. Яшина, акад. В.В. Козлов), научных семинарах ДФВУ г. Владивосток (рук. проф. Е.А. Нурминский), УФУ г. Екатеринбург (рук. проф. Н.Н. Субботина), БФУ г. Калининград (рук. проф. С.В. Мациевский), Иннополиса г. Иннополис (рук. доц. Я.А. Холодов). Результаты докладывались также на научном семинаре института К. Вейерштрасса в Берлине (рук. проф. В.Г. Спокойный).

Результаты диссертации докладывались на многочисленных международных конференциях и школах: Traffic and granular flows, Moscow 2011; OPTIMA (Optimization and applications) 2012–2015; ММРО (Математические методы распознавания образов) 2013, 2015; ИОИ (Интеллектуализация обработки информации) 2012, 2014; ИТиС (Информационные технологии и системы) 2013– 2015; ISMP (International Symposium on Mathematical Programming) 2012, 2015; International Conference Network Analysis (NET2015), Нижний Новгород, 2015; Летняя школа "Современная математика", Дубна 2011–2014, 2016; Сириус-2016; Традиционная математическая школа "Управление, информация и оптимизация" Б.Т. Поляка 2012–2016; Summer School on Operation Research and Application, Нижний Новгород, 2015; Байкальские чтения Алексея Савватеева, 2014–2016; Skoltech Deep Machine Intelligence workshop, 2016. На конференции "Moscow International Conference on Operation Research (ORM)", 2013, соискатель руководил работой транспортной секции.

По тематике диссертации соискателем организовано несколько научных семинаров: "Математическое моделирование транспортных потоков" (совместный семинар МФТИ–НМУ), "Стохастический анализ в задачах" (совместный семинар МФТИ–НМУ), Математический кружок (МФТИ), на которых в 2012–2016 гг. докладывались результаты работы.

**Благодарности**

Настоящая работа была инициирована в конце 2007 года проф. А.А. Шананиным, предложившим соискателю заняться математическим моделированием транспортных потоков и параллельно начать читать студентам МФТИ одноименный курс (такой курс читается с 2008 года). Полезным при написании работы оказался фундамент численных методов выпуклой оптимизации, заложенный во время обучения на базовой кафедре ФУПМ МФТИ в ВЦ РАН. Важным этапом для научного роста соискателя стало создание в 2011 г. на ФУПМ МФТИ лаборатории Структурных методов анализа данных в предсказательном моделировании (рук. проф. В.Г. Спокойный). Благодаря лаборатории и лично проф. В.Г. Спокойному и акад. А.П. Кулешову у соискателя появились новые возможности и новые научные контакты, которые и определили окончательный выбор темы диссертации (основные результаты этой работы были получены в МФТИ и ИППИ РАН в период апрель 2012 года – апрель 2016 года). Принципиально важную роль сыграло научное общение с

## Структура и объем работы

Диссертация состоит из введения, шести глав, заключения, одного приложения и списка литературы. Полный объем диссертации составляет 487 страниц, список литературы содержит 331 наименование.

## Основное содержание работы

Во **Введении** дается общая характеристика диссертации, обосновывается актуальность темы исследования и степень ее разработанности. Основные аргументы из Введения приведены в первой части автореферата.

**Глава 1** является самой большой главой диссертации. В этой главе выводятся вариационные принципы, описывающие равновесия в транспортных сетях. Ниже мы изложим основную конструкцию в упрощенном варианте (без эволюционной составляющей). Первый параграф этой главы содержит обзор основных результатов первой главы и диссертации (с точки зрения транспортных приложений) в целом.

Рассмотрим транспортную сеть, заданную ориентированным графом $\Gamma^1 = \langle V^1, E^1 \rangle$. Часть его вершин $O^1 \subseteq V^1$ является источниками, часть стоками $D^1 \subseteq V^1$. Множество пар источник-сток, обозначим $OD^1 \subseteq O^1 \otimes D^1$. Пусть каждой паре $w^1 \in OD^1$ соответствует своя корреспонденция: $d_{w^1}^1 := d_{w^1}^1 \cdot N$ ($N \gg 1$) пользователей, которые хотят в единицу времени перемещаться из источника в сток, соответствующих заданной корреспонденции $w^1$. Пусть ребра $\Gamma^1$ разделены на два типа $E^1 = \tilde{E}^1 \amalg \bar{E}^1$. Ребра типа $\tilde{E}^1$ характеризуются неубывающими функциями затрат $\tau_{e^1}^1\left(f_{e^1}^1\right) := \tau_{e^1}^1\left(f_{e^1}^1/N\right)$. Затраты $\tau_{e^1}^1\left(f_{e^1}^1\right)$ несут те пользователи, которые используют в своем пути ребро $e^1 \in \tilde{E}^1$, в предположении, что поток пользователей по этому ребру равен $f_{e^1}^1$. Пары вершин, задающие ребра типа $\bar{E}^1$, являются, в свою очередь, парами источник-сток $OD^2$ (с корреспонденциями $d_{w^2}^2 = f_{e^1}^1$, $w^2 = e^1 \in \bar{E}_1$) в транспортной сети следующего уровня, $\Gamma^2 = \langle V^2, E^2 \rangle$, ребра которой, в свою очередь, разделены на два типа $E^2 = \tilde{E}^2 \amalg \bar{E}^2$. Ребра типа $\tilde{E}^2$ характеризуются неубывающими функциями затрат $\tau_{e^2}^2\left(f_{e^2}^2\right) := \tau_{e^2}^2\left(f_{e^2}^2/N\right)$. Затраты $\tau_{e^2}^2\left(f_{e^2}^2\right)$ несут те пользователи, которые используют в своем пути ребро $e^2 \in \tilde{E}^2$, в предположении, что поток пользователей по этому ребру равен $f_{e^2}^2$. Пары вершин, задающие ребра типа $\bar{E}^2$, являются, в свою очередь, парами источник-сток $OD^3$ (с корреспонденциями $d_{w^3}^3 = f_{e^2}^2$, $w^3 = e^2 \in \bar{E}^2$) в транспортной сети более высокого уровня, $\Gamma^3 = \langle V^3, E^3 \rangle$, ... и т.д. Будем считать, что всего имеется $m$ уровней: $\tilde{E}^m = E^m$. Обычно в приложениях число $m$ небольшое: 1 – 4.

Каждый пользователь в графе $\Gamma^1$ выбирает путь $p_{w^1}^1 \in P_{w^1}^1$ (последовательный набор проходимых пользователем ребер), соответствующий его корреспонденции $w^1 \in OD^1$ ($P_{w^1}^1$ – множество всех путей, отвечающих в $\Gamma^1$ корреспонденции $w^1$). Задав $p_{w^1}^1$ можно однозначно восстановить ребра типа $\bar{E}^1$, входящие в этот путь. На каждом из этих ребер $w^2 \in \bar{E}^1$ пользователь может выбирать свой путь $p_{w^2}^2 \in P_{w^2}^2$ ($P_{w^2}^2$ – множество всех путей, отвечающих в $\Gamma^2$ корреспонденции $w^2$), ... и т.д. Пусть каждый пользователь сделал свой выбор. Обозначим через $x_{p^1}^1$ величину потока пользователей по пути $p^1 \in P^1 = \coprod\limits_{w^1 \in OD^1} P_{w^1}^1$, $x_{p^2}^2$ – величина потока пользователей по пути $p^2 \in P^2 = \coprod\limits_{w^2 \in OD^2} P_{w^2}^2$, ... и т.д. Заметим, что

$$x_{p_{w^k}^k}^k \geq 0,\ p_{w^k}^k \in P_{w^k}^k,\ \sum_{p_{w^k}^k \in P_{w^k}^k} x_{p_{w^k}^k}^k = d_{w^k}^k,\ w^k \in OD^k,\ k=1,...,m.$$

Отметим, что здесь и везде в дальнейшем

$$w^{k+1}\left(= e^k\right) \in OD^{k+1}\left(= \bar{E}^k\right),\ d_{w^{k+1}}^{k+1} = f_{e^k}^k,\ k=1,...,m-1.$$

Введем для графа $\Gamma^k$ и множества путей $P^k$ матрицу Кирхгофа

$$\Theta^k = \left\|\delta_{e^k p^k}\right\|_{e^k \in E^k, p^k \in P^k},\ \delta_{e^k p^k} = \begin{cases} 1, & e^k \in p^k \\ 0, & e^k \notin p^k \end{cases},\ k=1,...,m.$$

Тогда вектор потоков на ребрах $f^k$ на графе $\Gamma^k$ однозначно определяется вектором потоков на путях $x^k = \left\{x_{p^k}^k\right\}_{p^k \in P^k}$:

$$f^k = \Theta^k x^k,\ k=1,...,m.$$

Обозначим через

$$x = \left\{x^k\right\}_{k=1}^m,\ f = \left\{f^k\right\}_{k=1}^m,\ \Theta = \mathrm{diag}\left\{\Theta^k\right\}_{k=1}^m.$$

Введем $E = \coprod_{k=1}^m \tilde{E}^k$ и положим $t = \{t_e\}_{e \in E}$. Определим индуктивно по $k < m$ функции (база индукции: $g_{p^m}^m(t) = \sum_{e^m \in E^m} \delta_{e^m p^m} t_{e^m}$)

$$g_{p^k}^k(t) = \sum_{e^k \in \tilde{E}^k} \delta_{e^k p^k} t_{e^k} - \sum_{e^k \in \bar{E}^k} \delta_{e^k p^k} \gamma^{k+1} \psi_{e^k}^{k+1}\left(t/\gamma^{k+1}\right)$$

– "длина" пути $p^k$ в графе $\Gamma^k$, ребра $e^k \in \tilde{E}^k$ которого имеют веса $t_{e^k}$, а ребра $e^k \in \bar{E}^k$ веса $\gamma^{k+1} \psi_{e^k}^{k+1}\left(t/\gamma^{k+1}\right)$, где параметр $\gamma^{k+1} \geq 0$ характеризует ограниченную рациональность (масштаб шума) пользователей на уровне $k$, а

$$\psi_{e^k}^{k+1}(t) = \psi_{w^{k+1}}^{k+1}(t) = \ln\left(\sum_{p^{k+1} \in P_{w^{k+1}}^{k+1}} \exp\left(-g_{p^{k+1}}^{k+1}(t)\right)\right).$$

Предположим, что каждый пользователь $l$ транспортной сети, использующий корреспонденцию $w^k \in OD^k$ на уровне $k$ (ребро $e^{k-1}\left(=w^k\right) \in \bar{E}^{k-1}$ на уровне $k-1$), выбирает маршрут следования $p^k \in P_{w^k}^k$ на уровне $k$, если

$$p^k = \arg\max_{q^k \in P_{w^k}^k}\left\{-g_{q^k}^k(t) + \xi_{q^k}^{k,l}\right\},$$

где независимые случайные величины $\xi_{q^k}^{k,l}$ имеют одинаковое двойное экспоненциальное распределение, также называемое распределением Гумбеля:

$$P\left(\xi_{q^k}^{k,l} < \zeta\right) = \exp\left\{-e^{-\zeta/\gamma^k - E}\right\}.$$

Такое распределение возникает в данном контексте, например, если при принятии решения водитель/пользователь собирает информацию с большого числа разных (независимых) зашумленных источников, ориентируясь на худшие прогнозы по каждому из путей. Отметим также, что если взять $E \approx 0.5772$ – константа Эйлера, то

$$E\left[\xi_{q^k}^{k,l}\right] = 0, \ D\left[\xi_{q^k}^{k,l}\right] = \left(\gamma^k\right)^2 \pi^2 / 6.$$

Распределение Гиббса (логит-распределение) пользователей по путям

$$x_{p^k}^k = d_{w^k}^k \frac{\exp\left(-g_{p^k}^k(t)/\gamma^k\right)}{\sum_{\tilde{p}^k \in P_{w^k}^k} \exp\left(-g_{\tilde{p}^k}^k(t)/\gamma^k\right)}, \ p^k \in P_{w^k}^k, \ w^k \in OD^k, \ k = 1,...,m. \quad (1)$$

получается в пределе (см., например, W. Sandholm, 2010), когда число агентов на каждой корреспонденции $w^k \in OD^k$, $k = 1,...,m$ стремится к бесконечности, т.е. $N \to \infty$ (случайность исчезает и описание переходит на средние величины). Полезно также в этой связи иметь в виду, что

$$\gamma^k \psi_{w^k}^k\left(t/\gamma^k\right) = E_{\left\{\xi_{p^k}^{k,l}\right\}_{p^k \in P_{w^k}^k}} \left[\max_{p^k \in P_{w^k}^k} \left\{-g_{p^k}^k(t) + \xi_{p^k}^{k,l}\right\}\right],$$

$$f = \Theta x = -\nabla \psi^1\left(t/\gamma^1\right), \ \psi^1(t) = \sum_{w^1 \in OD^1} d_{w^1}^1 \psi_{w^1}^1(t). \quad (1')$$

Если каждый пользователь сориентирован на вектор затрат $t$ на ребрах $E$, одинаковый для всех пользователей, и на каждом уровне принятия решения пытается выбрать кратчайший путь, исходя из зашумленной информации и исходя из усреднения деталей более высоких уровней, то такое поведение пользователей в пределе, когда их число стремится к бесконечности $N \to \infty$, приводит к описанию распределения пользователей по путям/ребрам (1). Такое усреднение можно обосновать, если, например, ввести разный масштаб времени (частот принятия решений) на разных уровнях. Можно просто постулировать, что пользователь действует так, как это принято в моделях типа Nested Logit (Y. Sheffi, 1985). *Равновесная конфигурация характеризуется тем, что по вектору t вычисляется согласно формуле (1) такой вектор* $f = \Theta x$ *(см. также (1')), что имеет место соотношение* $t = \left\{\tau_e(f_e)\right\}_{e \in E}$.

Введя $\sigma_e(f_e) = \int_0^{f_e} \tau_e(z)dz$, $\sigma_e^*(t_e) = \max_{f_e}\{f_e t_e - \sigma_e(f_e)\}$, получим

$$\frac{d\sigma_e^*(t_e)}{dt_e} = \frac{d}{dt_e}\max_{f_e}\left\{f_e t_e - \int_0^{f_e} \tau_e(z)dz\right\} = f_e: t_e = \tau_e(f_e),\ e \in E.$$

Это наблюдение позволяет установить следующий результат.

**Теорема 1 (вариационный принцип).** *Поиск равновесных $x$, $f$, $t$ (неподвижной точки) приводит к следующей паре задач (далее используется обозначение $\mathrm{dom}\,\sigma_e^*$ – область определения сопряженной к $\sigma_e^*$ функции)*

$$\min_{f,x}\{\Psi(x,f):\ f = \Theta x,\ x \in X\} =$$
$$= -\min_{t \in \{\mathrm{dom}\,\sigma_e^*\}_{e \in E}} \left\{\gamma^1 \psi^1(t/\gamma^1) + \sum_{e \in E} \sigma_e^*(t_e)\right\}, \quad (2)$$

*где*

$$\Psi(x,f) := \Psi^1(x) = \sum_{e^1 \in \tilde{E}^1} \sigma_{e^1}^1(f_{e^1}^1) + \Psi^2(x) + \gamma^1 \sum_{w^1 \in OD^1} \sum_{p^1 \in P_{w^1}^1} x_{p^1}^1 \ln(x_{p^1}^1/d_{w^1}^1),$$

$$\Psi^2(x) = \sum_{e^2 \in \tilde{E}^2} \sigma_{e^2}^2(f_{e^2}^2) + \Psi^3(x) + \gamma^2 \sum_{w^2 \in \bar{E}^1} \sum_{p^2 \in P_{w^2}^2} x_{p^2}^2 \ln(x_{p^2}^2/d_{w^2}^2),\ d_{w^2}^2 = f_{w^2}^1,$$

..................................................................................................

$$\Psi^k(x) = \sum_{e^k \in \tilde{E}^k} \sigma_{e^k}^k(f_{e^k}^k) + \Psi^{k+1}(x) + \gamma^k \sum_{w^k \in \bar{E}^{k-1}} \sum_{p^k \in P_{w^k}^k} x_{p^k}^k \ln(x_{p^k}^k/d_{w^k}^k),\ d_{w^{k+1}}^{k+1} = f_{w^{k+1}}^k,$$

..................................................................................................

$$\Psi^m(x) = \sum_{e^m \in E^m} \sigma_{e^m}^m(f_{e^m}^m) + \gamma^m \sum_{w^m \in \bar{E}^{m-1}} \sum_{p^m \in P_{w^m}^m} x_{p^m}^m \ln(x_{p^m}^m/d_{w^m}^m),\ d_{w^m}^m = f_{w^m}^{m-1}.$$

Слагаемое $\gamma^1 \psi^1(t/\gamma^1)$ в двойственной задаче (2) имеет равномерно ограниченную константу Липшица градиента в 2-норме:

$$L_2 \leq \frac{1}{\min_{k=1,\ldots,m} \gamma^k} \sum_{w^1 \in OD^1} d_{w^1}^1 \cdot (l_{w^1})^2,$$

где $l_{w^1}$ – число рёбер (всех уровней) в самом длинном пути из корреспонденции $w^1$. Эта гладкость теряется при $\gamma^k \to 0+$ (если $\gamma^k \to 0+$, то на уровне $k$ пользователи рациональны, и соответствующий блок отвечает равновесию Нэша, а не стохастическому равновесию): $-\lim_{\gamma^k \to 0+} \gamma^k \psi_{w^k}^k(t/\gamma^k) = \min_{p^k \in P_{w^k}^k} g_{p^k}^k(t)$ – длина кратчайшего пути на корреспонденции $w^k \in OD^k$ в графе $\Gamma^k$.

Во многих приложениях выбирают (BPR-функции)

$$\tau_e(f_e) := \tau_e^\mu(f_e) = \overline{t}_e \cdot \left(1 + \gamma \cdot \left(\frac{f_e}{\overline{f}_e}\right)^{\frac{1}{\mu}}\right),$$

с $\mu = 0.25$, где $\overline{t}_e$ – затраты на прохождение ребра $\overline{t}$, когда ребро свободно, $\overline{f}_e$ – пропускная способность ребра. В пределе Стабильной Динамики

$$\tau_e^\mu(f_e) \xrightarrow[\mu \to 0+]{} \begin{cases} \overline{t}_e, & f_e < \overline{f}_e \\ \infty, & f_e > \overline{f}_e \end{cases},$$

получаем (в качестве $\tau_e^\mu(f_e)$ можно использовать произвольные гладкие монотонные функции с таким свойством, т.е. не обязательно BPR типа)

$$\sigma_e^*(t_e) = \lim_{\mu \to 0+} \max_{f_e} \left\{ f_e t_e - \int_0^{f_e} \tau_e^\mu(z) dz \right\} = \begin{cases} \overline{f}_e \cdot (t_e - \overline{t}_e), & t_e \geq \overline{t}_e \\ \infty, & t_e < \overline{t}_e \end{cases}.$$

При этом $\overline{f}_e - f_e$ в точности совпадает с множителем Лагранжа к возникающему в таком пределе дополнительному ограничению $t_e \geq \overline{t}_e$.

Все известные нам модели равновесного распределения транспортных потоков (в том числе многостадийные) могут быть получены из формулы (2) с помощью выбора $m$ и осуществления/не осуществления (в различных сочетаниях) описанных выше предельных переходов. В частности, выбирая $m = 1$ и вводя потенциалы притяжения/отталкивания районов (потенциалы Канторовича–Гавурина), из (2) можно получить энтропийную модель расчета матрицы корреспонденций. Выбирая $m = 1$ и полагая $\gamma^1 = \gamma \to 0+$, получим модель Бэкмана. Если, дополнительно, по всем ребрам осуществить предельный переход $\mu \to 0+$, то получим модель Стабильной Динамики, а если такой предельный переход осуществить только по части ребер, то получим Смешанную Модель, встречающуюся, например, при изучении тарифной политики грузоперевозок РЖД (А.А. Шананин и др., 2014). Выбирая $m = 2$, $\gamma^2 \to 0+$, получим обычную (обычную не по форме записи, а по конечному результату – описанию равновесия) двухстадийную транспортную модель.

В первой главе также намечаются различные подходы к численному решению задачи (2) в общем случае и в различных специальных случаях. Более подробно эти подходы рассматриваются в последующих главах.

В заключение обзора результатов первой главы отметим, что для решения задачи (2) численным методом (например, универсальным методом тре-угольника из главы 2) на каждой итерации метода необходимо вычислять $\nabla \psi^1(t/\gamma^1)$, а для ряда методов и $\psi^1(t/\gamma^1)$ (например, для адаптивных методов, настраивающихся на параметры гладкости задачи). Используя быстрое автоматическое дифференцирование, можно показать (с помощью сглаженного варианта метода Форда–Беллмана, предложенного Ю.Е. Нестеровым, 2007; при $\gamma \to 0+$ стоит использовать обычный метод Форда–Беллмана или алгоритм Дейкстры), что для этого достаточно сделать $\mathrm{O}\left(\left|O^1\right|\left|E\right|\max_{w^1 \in OD^1} l_{w^1}\right)$ арифметических операций.

В **Главе 2** рассматриваются основные (оптимальные по Немировскому–Юдину, 1979) численные методы решения (гладких) задач выпуклой (композитной стохастической) оптимизации на множествах простой структуры в случае наличия шума не случайной природы (также говорят – в случае неточного оракула). Первый параграф этой главы содержит обзор современного состояния стохастических градиентных методов с неточным оракулом. В этот параграф, в том числе, вошли результаты, которые более подробно излагаются в последующих частях диссертации. Основным результатом этой главы является построение на базе варианта БГМ (Метода Треугольника, Ю.Е. Нестеров, 2016) универсального метода, способного работать с сильно выпуклыми постановками задач и со стохастическими постановками. Здесь и далее под выпуклой задачей понимается задача оптимизации, в которой все функции и множества выпуклые. Сильно выпуклая постановка дополнительно предполагает сильную выпуклость функционала задачи.

Рассматривается задача (включающая в себя двойственную задачу (2)) выпуклой композитной оптимизации (размерность вектора $x$ велика)

$$F(x) = f(x) + h(x) \to \min_{x \in Q}. \qquad (3)$$

Положим $R^2 = V(x_*, y^0)$, где прокс-расстояние (расстояние Брэгмана) определяется формулой $V(x, z) = d(x) - d(z) - \langle \nabla d(z), x - z \rangle$; прокс-функция $d(x) \geq 0$ ($d(y^0) = 0$, $\nabla d(y^0) = 0$) считается сильно выпуклой относительно выбранной нормы $\|\ \|$ ($\|\ \|_*$ – двойственная норма), с константой сильной выпуклости $\geq 1$; $x_*$ – решение задачи (3) (если решение не единственно, то выбирается то, которое доставляет минимум $V(x_*, y^0)$), $\varepsilon$ – желаемая точность.

**Предположение 1.** *Пусть*

$$\|\nabla f(y) - \nabla f(x)\|_* \leq L_\nu \|y-x\|^\nu, \ \nu \in [0,1].$$

**Предположение 2.** *Пусть $f(x)$ – $\mu$-сильно выпуклая функция в норме $\|\ \|$, т.е. для любых $x, y \in Q$ имеет место неравенство*

$$f(y) + \langle \nabla f(y), x-y \rangle + \frac{\mu}{2}\|x-y\|^2 \leq f(x).$$

Введем $\tilde{\omega}_n = \sup_{x,y \in Q} 2V(x,y)/\|y-x\|^2$. В евклидовом случае выбирают $\|\ \| = \|\ \|_2$, $V(x,y) = \|x-y\|_2^2/2$, что дает $\tilde{\omega}_n = 1$. Положим ($\tilde{\mu} = \mu/\tilde{\omega}_n$)

$$\varphi_0(x) = V(x, y^0) + \alpha_0\left[f(y^0) + \langle \nabla f(y^0), x-y^0 \rangle + \tilde{\mu}V(x, y^0) + h(x)\right],$$

$$\varphi_{k+1}(x) = \varphi_k(x) + \alpha_{k+1}\left[f(y^{k+1}) + \langle \nabla f(y^{k+1}), x-y^{k+1} \rangle + \tilde{\mu}V(x, y^k) + h(x)\right].$$

Положим $A_0 = \alpha_0 = 1/L_0^0 = 1$, $k = 0$, $j_0 = 0$. До тех пор пока

$$f(y^0) + \langle \nabla f(y^0), x^0 - y^0 \rangle + \frac{L_0^{j_0}}{2}\|x^0 - y^0\|^2 + \frac{\alpha_0}{2A_0}\varepsilon < f(x^0),$$

где $x^0 := u^0 := \arg\min_{x \in Q} \varphi_0(x)$, $(A_0 :=)\alpha_0 := \dfrac{1}{L_0^{j_0}}$, выполнять $j_0 := j_0 + 1$; $L_0^{j_0} := 2^{j_0}L_0^0$.

## Универсальный Метод Треугольника (УМТ)

1. $L_{k+1}^0 = L_k^{j_k}/2$, $j_{k+1} = 0$.

2. $\begin{cases} \alpha_{k+1} := \dfrac{1 + A_k\tilde{\mu}}{2L_{k+1}^{j_{k+1}}} + \sqrt{\dfrac{1 + A_k\tilde{\mu}}{4\left(L_{k+1}^{j_{k+1}}\right)^2} + \dfrac{A_k \cdot (1 + A_k\tilde{\mu})}{L_{k+1}^{j_{k+1}}}}, A_{k+1} := A_k + \alpha_{k+1}; \\ y^{k+1} := \dfrac{\alpha_{k+1}u^k + A_k x^k}{A_{k+1}}, u^{k+1} := \arg\min_{x \in Q} \varphi_{k+1}(x), x^{k+1} := \dfrac{\alpha_{k+1}u^{k+1} + A_k x^k}{A_{k+1}}. \end{cases}$ (*)

До тех пор пока

$$f(y^{k+1}) + \langle \nabla f(y^{k+1}), x^{k+1} - y^{k+1} \rangle + \frac{L_{k+1}^{j_{k+1}}}{2}\|x^{k+1} - y^{k+1}\|^2 + \frac{\alpha_{k+1}}{2A_{k+1}}\varepsilon < f(x^{k+1}),$$

выполнять

$$j_{k+1} := j_{k+1} + 1;\ L_{k+1}^{j_{k+1}} = 2^{j_{k+1}} L_{k+1}^0;\ (*).$$

3. Если не выполнен критерий останова, то $k := k+1$ и **go to** 1.

**Теорема 2.** *Пусть предположение 1 выполняется хотя бы для $v = 0$, и справедливо предположение 2 с $\mu \geq 0$ (допускается брать $\mu = 0$). Тогда УМТ для задачи (3) сходится согласно оценке*

$$F(x^N) - \min_{x \in Q} F(x) \leq \varepsilon, \text{ для любого } N \geq N(\varepsilon), \text{ где}$$

$$N(\varepsilon) = \min \left\{ \inf_{v \in [0,1]} \left( \frac{L_v \cdot (16R)^{1+v}}{\varepsilon} \right)^{\frac{2}{1+3v}}, \inf_{v \in [0,1]} \left\{ \left( \frac{8 L_v^{\frac{2}{1+v}} \tilde{\omega}_n}{\mu \varepsilon^{\frac{1-v}{1+v}}} \right)^{\frac{1+v}{1+3v}} \ln^{\frac{2+2v}{1+3v}} \left( \frac{16 L_v^{\frac{4+6v}{1+v}} R^2}{(\mu/\tilde{\omega}_n)^{\frac{1+v}{1+3v}} \varepsilon^{\frac{5+7v}{2+6v}}} \right) \right\} \right\}. \quad (4)$$

*При этом среднее число вычислений значения функции на одной итерации будет равно* 4*, а градиента функции* – 2.

Заметим, что при этом для любого $k = 0, 1, 2, ...$ имеют место неравенства

$$\|u^k - x_*\|^2 \leq 2R^2,\ \max\left\{\|x^k - x_*\|^2, \|y^k - x_*\|^2\right\} \leq 4R^2 + 2\|x^0 - y^0\|^2.$$

Если inf достигается при $v = 0$, то УМТ соответствует (с точностью до логарифмического множителя) по оценке скорости сходимости МЗС, а если при $v = 1$, то БГМ.

Предположим теперь, что вместо настоящих градиентов доступны только стохастические градиенты $\nabla f(x) \to \nabla f(x, \xi)$ (можно обобщить на случай, когда вместо значений функции доступны только реализации значений функции $f(x) \to f(x, \xi)$).

**Предположение 3.** *Пусть для всех $x \in Q$*

$$E_\xi [\nabla f(x, \xi)] = \nabla f(x) \text{ и } E_\xi \left[ \|\nabla f(x, \xi) - \nabla f(x)\|_*^2 \right] \leq D.$$

Введём обозначение (mini-batch'инг)

$$\bar{\nabla}^m f(x) = \frac{1}{m} \sum_{k=1}^m \nabla f(x, \xi^k),$$

где $\xi^k$ – независимые одинаково распределенные (так же как $\xi$) случайные величины. Переопределим последовательность $\varphi_{k+1}(x)$, заменив $\nabla f(y^{k+1})$ на $\overline{\nabla}^{m_{k+1}} f(y^{k+1})$ (аналогично с $\varphi_0(x)$). Тогда если дополнительно в условиях теоремы 2 имеет место предположение 3 и на шаге 2 УМТ ввести $m_{k+1} := \left\lceil 8DA_{k+1} / L_{k+1}^{j_{k+1}} \alpha_{k+1} \varepsilon \right\rceil$, заменив условие выхода из цикла следующим

$$f(y^{k+1}) + \left\langle \overline{\nabla}^{m_{k+1}} f(y^{k+1}), x^{k+1} - y^{k+1} \right\rangle + \frac{L_{k+1}^{j_{k+1}}}{2} \left\| x^{k+1} - y^{k+1} \right\|^2 + \frac{\alpha_{k+1}}{2A_{k+1}} \varepsilon < f(x^{k+1}),$$

то оценка (4) видоизменится следующим образом: $N(\varepsilon) \to 2N(\varepsilon/4)$. При этом среднее число вычислений значения функции на одной итерации по-прежнему будет равно $4$. А оценка числа обращений за стохастическим градиентом для достижения точности по функции (в среднем) $\varepsilon$ примет вид (с точностью до множителя $\sim \ln n$ в случае неевклидовой прокс-структуры)

$$2 \cdot \min \left\{ \frac{64DR^2}{\varepsilon^2}, \frac{8D\tilde{\omega}_n}{\mu \varepsilon} \ln \left( \frac{8L_0^{j_0} R^2}{\varepsilon} \right) \right\} + 4N(\varepsilon/4). \qquad (5)$$

Оценки (4), (5) сохранят свой вид (немного увеличатся числовые коэффициенты), если считать, что вместо оракула, выдающего градиент (несмещенную оценку градиента) и значение функции, у нас есть доступ только к введенному Деволдером–Нестеровым–Глинёром $(\delta, L, \mu)$-оракулу с $\delta = \mathrm{O}(\varepsilon/N(\varepsilon))$ и $L = \mathrm{O}\left( \max_{k=0,...,N} L_k^{j_k} \right)$. Этот оракул на запрос, в котором указывается только одна точка $x$, выдает такую пару $(f_\delta(x), g_\delta(x, \xi))$ (можно обобщить на случай, когда $f_\delta(x) \to f_\delta(x, \xi)$), что для всех $x \in Q$

$$E_\xi \left[ \left\| g_\delta(x, \xi) - E_\xi [g_\delta(x, \xi)] \right\|_*^2 \right] \le D,$$

и для любых $x, y \in Q$

$$\frac{\mu}{2} \|y - x\|^2 \le f(y) - f_\delta(x) - \left\langle E_\xi [g_\delta(x, \xi)], y - x \right\rangle \le \frac{L}{2} \|y - x\|^2 + \delta.$$

С точностью до числовых и логарифмических множителей оценки (4), (5) оптимальны, т.е. без дополнительных уточнений постановки задачи эти оценки не могут быть улучшены на соответствующих классах (гладкости /

сильной выпуклости) функций $f(x)$ (Немировский–Юдин, 1979). Оптимальна и оценка $\delta = \mathrm{O}(\varepsilon/N(\varepsilon))$ на уровень допустимого шума (Деволдер–Нестеров–Глинёр, 2013).

Рассмотрим конкретные примеры приложений, описанных в теореме 2 (и последующем тексте) результатов. В этих примерах известно, что существует $L = L_1 < \infty$ (см. обозначения в предположении 1).

**Пример 1 (min max-задача).** Рассматривается задача поиска седловой точки (говорят также, что функционал представим в форме Лежандра)

$$f(x) = \max_{\|y\|_2 \leq R_y}\{G(y) + \langle By, x\rangle\} \to \min_{\|x\|_2 \leq R_x},$$

где функция $G(y)$ – сильно вогнутая с константой $\mu$ относительно 2-нормы и константой Липшица градиента $L_G$ (также в 2-норме). Тогда функция $f(x)$ будет гладкой, с константой Липшица градиента в 2-норме $L_f = \sigma_{\max}(B)/\mu$. Казалось бы, что можно решить задачу минимизации функции $f(x)$ за $\mathrm{O}\left(\sqrt{\sigma_{\max}(B)R_x^2/(\mu\varepsilon)}\right)$ итераций, где $\varepsilon$ – желаемая точность по функции $f(x)$. Но это возможно только, если мы можем абсолютно точно находить $\nabla f(x) = By^*(x)$, где $y^*(x)$ – решение вспомогательной задачи максимизации по $y$ при заданном $x$. В действительности, мы можем решать эту задачу (при различных $x$) лишь приближенно. Если мы решаем вспомогательную задачу УМТ с точностью $\delta/2$ (на это потребуется $\mathrm{O}\left(\sqrt{L_G/\mu}\ln\left(L_G R_y^2/\delta\right)\right)$ итераций), то пара $\left(G(y_{\delta/2}(x)) + \langle By_{\delta/2}(x), x\rangle, By_{\delta/2}(x)\right)$, где $y_{\delta/2}(x)$ – $\delta/2$-решение вспомогательной задачи, будет $(\delta, 2L_f, 0)$-оракулом (Деволдер–Нестеров–Глинёр, 2013). Выбирая $\delta = \mathrm{O}\left(\varepsilon\sqrt{\varepsilon/(L_f R_x^2)}\right)$, получим после

$$\mathrm{O}\left(\sqrt{\frac{L_G \sigma_{\max}(B) R_x^2}{\mu^2 \varepsilon}} \ln\left(\frac{L_f L_G R_x^2 R_y^2}{\varepsilon}\right)\right)$$

итераций (на итерациях производится умножение матрицы $B$ на вектор/строчку и вычисление градиента $G(y)$) $\varepsilon$-решение задачи минимизации $f(x)$. Отметим, что если не использовать сильную вогнутость функции $G(y)$, то для получения пары $(x^N, y^N)$, удовлетворяющей неравенству (что, фактически, отвечает решению исходной задачи с точностью $\varepsilon$)

$$\max_{\|y\|_2 \leq R_y} \{G(y) + \langle By, x^N \rangle\} - \min_{\|x\|_2 \leq R_x} \{G(y^N) + \langle By^N, x \rangle\} \leq \varepsilon,$$

потребуется $\Omega\left(\max\{L_G R_y^2, \sigma_{\max}(B) R_x R_y\}/\varepsilon\right)$ итераций.

**Пример 2 (min min-задача).** Пусть $f(x) = \min_{y \in Q} \Phi(y, x)$, где $Q$ — ограниченное выпуклое множество, а $\Phi(y, x)$ — такая достаточно гладкая, выпуклая по совокупности переменных функция, что при $y, y' \in Q$, $x, x' \in \mathbb{R}^n$

$$\|\nabla \Phi(y', x') - \nabla \Phi(y, x)\|_2 \leq L \|(y', x') - (y, x)\|_2.$$

Пусть для произвольного $x$ можно найти такой $\tilde{y}(x) \in Q$, что

$$\max_{z \in Q} \langle \nabla_y \Phi(\tilde{y}(x), x), \tilde{y}(x) - z \rangle \leq \delta.$$

Тогда для любых $x, x' \in \mathbb{R}^n$

$$\Phi(\tilde{y}(x), x) - f(x) \leq \delta, \quad \|\nabla f(x') - \nabla f(x)\|_2 \leq L \|x' - x\|_2,$$

и $\left(\Phi(\tilde{y}(x), x) - 2\delta, \nabla_y \Phi(\tilde{y}(x), x)\right)$ будет $(6\delta, 2L, 0)$-оракулом для $f(x)$.

Этот пример оказался полезным для решения двойственной задачи (2) с помощью УМТ (с $\mu = 0$).

**Пример 3 (восстановление матрицы корреспонденций в компьютерной сети по замерам потоков на линках / ребрах). a)** Рассмотрим следующую задачу выпуклой композитной оптимизации (вместо ограничения $\sum_{k=1}^n x_k = 1$ можно рассматривать ограничение $\sum_{k=1}^n x_k \leq 1$):

$$F(x) = \frac{1}{2}\|Ax - b\|_2^2 + \mu \sum_{k=1}^n x_k \ln x_k \to \min_{\sum_{k=1}^n x_k = 1,\, x \geq 0}.$$

Разберем два случая: а) $0 < \mu \ll \varepsilon/(2\ln n)$ – мало (сильную выпуклость композита в 1-норме можно не учитывать); б) $\mu \gg \varepsilon/(2\ln n)$ – достаточно большое (сильную выпуклость композита в 1-норме необходимо учитывать). Выберем норму в прямом пространстве $\|\ \| = \|\ \|_1$. Положим

$$f(x) = \frac{1}{2}\|Ax - b\|_2^2,\ h(x) = \mu\sum_{k=1}^{n} x_k \ln x_k,\ Q = S_n(1) = \left\{x \geq 0 : \sum_{k=1}^{n} x_k = 1\right\},$$

$L = \max_{k=1,\ldots,n}\|A^{\langle k\rangle}\|_2^2$, где $A^{\langle k\rangle}$ – $k$-й столбец матрицы $A$. Для случая а) можно выбирать $d(x) = \ln n + \sum_{k=1}^{n} x_k \ln x_k$. Тогда $V(x, z) = \sum_{k=1}^{n} x_k \ln(x_k/z_k)$, $R^2 \leq \ln n$. Здесь имеет место ситуация, когда композит совпадает по форме с прокс-расстоянием (энтропийного типа), и шаг итерации УМТ (с $\mu = 0$) осуществим по явным формулам. Таким образом, стоимость итерации (число арифметических операций типа умножения двух *double* чисел на итерации) будет $\mathrm{O}(nnz(A))$, где $nnz(A)$ – число ненулевых элементов в матрице $A$ (считаем, что это число $\geq n$). Оценка числа итераций будет иметь вид (см. теорему 2)

$$N = \mathrm{O}\left(\sqrt{\frac{\max_{k=1,\ldots,n}\|A^{\langle k\rangle}\|_2^2 \ln n}{\varepsilon}}\right).$$

К сожалению, использовать УМТ в случае б) нельзя, поскольку сильно выпуклым является композит $h(x)$, а не гладкая часть функционала $f(x)$. Однако можно на базе УМТ (с $\mu = 0$) построить с помощью рестартов оптимальный метод. Описываемые далее (во многом известные) конструкции позволяют переносить алгоритмы с сильно выпуклых задач на не сильно выпуклые задачи и наоборот, сохраняя при этом оптимальность. Эти конструкции (наряду с mini-batch'ингом и переходом к двойственной задаче) являются основными операциями, которые активно используются на протяжении всей диссертационной работы.

Введем семейство $\mu$-сильно выпуклых в норме $\|\ \|$ задач ($\mu > 0$)

$$F^\mu(x) = F(x) + \mu V(x, y^0) \to \min_{x \in Q}. \tag{6}$$

Пусть $F_*^\mu$ – оптимальное значение в задаче (6), а $F_*$ – в задаче (3).

**Утверждение 1 (регуляризация).** *Пусть*

$$\mu \leq \frac{\varepsilon}{2V(x_*, y^0)} = \frac{\varepsilon}{2R^2},$$

*и удалось найти $\varepsilon/2$-решение задачи (6), т.е. нашелся такой $x^N \in Q$, что*

$$F^\mu(x^N) - F_*^\mu \leq \varepsilon/2.$$

*Тогда*

$$F(x^N) - F_* \leq \varepsilon.$$

**Утверждение 2 (рестарты).** *Пусть справедливо предположение 1 с $\nu = 1$ ($L = L_1$), функция $F(x)$ – $\mu$-сильно выпуклая в норме $\|\ \|$. Пусть точка $x^{\bar{N}}(y^0)$ выдается УМТ (с $\mu = 0$), стартующим с точки $y^0$, после*

$$\bar{N} = \sqrt{\frac{16 L \omega_n}{\mu}}$$

*итераций, где (следует сравнить с введенным ранее $\tilde{\omega}_n$)*

$$\omega_n = \sup_{x \in Q} \frac{2V(x, y^0)}{\|x - y^0\|^2}.$$

*Положим $\left[x^{\bar{N}}(y^0)\right]^1 = x^{\bar{N}}(y^0)$, и определим по индукции*

$$\left[x^{\bar{N}}(y^0)\right]^{k+1} = x^{\bar{N}}\left(\left[x^{\bar{N}}(y^0)\right]^k\right),\ k = 1, 2, \ldots.$$

*При этом на $k+1$ перезапуске (рестарте) также корректируется прокс-функция (считаем, что так определенная функция корректно определена на $Q$ с сохранением свойства сильной выпуклости)*

$$d^{k+1}(x) = d\left(x - \left[x^{\bar{N}}(y^0)\right]^k + y^0\right) \geq 0,$$

*чтобы $d^{k+1}\left(\left[x^{\bar{N}}(y^0)\right]^k\right) = 0$, $\nabla d^{k+1}\left(\left[x^{\bar{N}}(y^0)\right]^k\right) = 0$. Тогда*

$$F\left(\left[x^{\bar{N}}(y^0)\right]^k\right) - F_* \leq \frac{\mu \|y^0 - x_*\|^2}{2^{k+1}}.$$

**Пример 3. б)** В этом случае следует использовать рестарт-технику (см. утверждение 2). Но для выбранной в п. а) функции $V(x, z)$ (расстояние Кульбака–Лейблера) процедура рестартов некорректна. Однако существует другой способ выбора прокс-функции

$$d(x) = \frac{1}{2(a-1)} \|x\|_a^2, \ a = \frac{2\ln n}{2\ln n - 1}.$$

В этом случае имеем $R^2 = \mathrm{O}(\ln n)$, $\omega_n = \mathrm{O}(\ln n)$. Сложность выполнения одной итерации (дополнительная к вычислению градиента гладкой части функционала $\mathrm{O}(nnz(A))$) определяется тем, насколько эффективно можно решить задачу следующего вида

$$\tilde{F}(x) = \langle c, x \rangle + \|x\|_a^2 + \bar{\mu} \sum_{k=1}^{n} x_k \ln x_k \to \min_{x \in S_n(1)}. \qquad (7)$$

В конце второй главы показывается (на базе концепции неточного оракула) с помощью перехода к двойственной задаче и ее (приближенного) решения с помощью прямо-двойственной версии метода эллипсоидов, что задачу (7) можно решить за $\mathrm{O}(n \ln^2(n/\varepsilon))$ арифметических операций, что в типичных ситуациях много меньше $\mathrm{O}(nnz(A))$. При этом (см. утверждение 2 и теорему 2) необходимое число таких итераций можно оценить следующим образом

$$N = \mathrm{O}\left( \sqrt{\frac{\max\limits_{k=1,\ldots,n} \|A^{\langle k \rangle}\|_2^2 \ln n}{\mu}} \left\lceil \log_2\left(\frac{\mu}{\varepsilon}\right) \right\rceil \right).$$

Заметим, что из этой формулы с помощью утверждения 1 можно получить (с точностью до $\sim \sqrt{\ln n}$) оценку примера 3 а).

Этот пример имел одной из своих целей продемонстрировать, что для большого класса задач отсутствие явных формул для шага итерации – не есть сколько-нибудь сдерживающие обстоятельство для использования метода.

В **Главе 3** исследуются различные способы восстановления приближенного решения прямой (двойственной), задачи, исходя из приближенного решения сопряженной ей двойственной (прямой) задачи. Необходимость такого исследования мотивирована задачей (2) и ее различными частными случаями. Большинство способов решения этой задачи исходят из ее двойственного представления. Таким образом, численный метод выдает только вектор $t = \{t_e\}$. Однако интересен также и вектор $f = \{f_e\}$ (а в ряде случаев и $x = \{x_p\}$). Кажется, что здесь можно воспользоваться формулой $t_e = \tau_e(f_e)$, для восстановления $f$. Однако при таком "очевидном" способе восстановления теряется контроль точности. Поэтому нужны способы с хорошими теоре-тическими гарантиями, изучению которых и посвящена данная глава. Далее мы опишем лишь одну из четырех "сюжетных линий", рассматриваемых в данной главе.

Пусть требуется решать задачу (включающую в себя прямую задачу (2))

$$g(x) \to \min_{Ax=b,\, x\in Q}, \qquad (8)$$

где функция $g(x)$ – 1-сильно выпуклая в $p$-норме $(1 \le p \le 2)$. Построим двойственную задачу (см. также пример 1)

$$f(y) = \max_{x\in Q}\{\langle y, b - Ax\rangle - g(x)\} \to \min_{y}. \qquad (9)$$

Во многих важных приложениях основной вклад в вычислительную сложность внутренней задачи максимизации дает умножение $Ax$ ($A^T y$). Это так, например, для сепарабельных функционалов

$$g(x) = \sum_{k=1}^{n} g_k(x_k)$$

и параллелепипедных ограничениях $Q$. В частности, это имеет место для задач энтропийно-линейного программирования (ЭЛП), возникающих при расчете матрицы корреспонденций, в которых явно выписана зависимость $x(y)$.

В общем случае внутренняя задача максимизации не решается точно (по явным формулам). Тем не менее, за счет сильной выпуклости $g(x)$ (аналогичное можно сказать в случае сепарабельности $g(x)$, но отсутствии сильной выпуклости) оценка точности решения этой вспомогательной задачи (на каждой итерации внешнего метода) учитывается в оценке сложности ее решения в виде логарифмической зависимости. Как следствие, оговорки о неточности оракула, выдающего градиент для внешней задачи минимизации $f(y)$, можно опустить. Аккуратный учет этих обстоятельств приводит лишь к логарифмическим поправкам в итоговых оценках сложности метода (см. также примеры 1, 3). Поэтому для большей наглядности мы далее в рассуждениях будем просто считать, что есть явная формула $x(y)$.

Положим $L = \max_{\|x\|_p \leq 1} \|Ax\|_2^2$. В частности, для задачи ЭЛП $p = 1$, т.е. $L = \max_{k=1,\ldots,n} \|A^{\langle k \rangle}\|_2^2$, где $A^{\langle k \rangle}$ — $k$-й столбец матрицы $A^{\langle k \rangle}$. Для задачи поиска вектора PageRank (см. главу 4 ниже) $p = 2$, т.е. $L = \lambda_{\max}(A^T A) = \sigma_{\max}(A)$.

Пусть УМТ (с $\mu = 0$) с $\|\ \| = \|\ \|_2$, $d(y) = \frac{1}{2}\|y\|_2^2$, $y^0 = 0$, для задачи (9) генерирует "модельные" точки $\{y^k\}$ (на основе которых строятся, по Ю.Е. Нестерову, модели функционала задачи $\varphi_k(y)$; в главе 2 эти точки также обозначали $\{y^k\}$), а выдает в итоге $\tilde{y}^N$ (соответствует $x^N$ главы 2). Положим

$$x^N = \sum_{k=0}^{N} \lambda_k x(y^k),\ \lambda_k = \alpha_k / A_N.$$

Поскольку ($x_*$ – решение задачи (8))

$$g(x^N) - g(x_*) \leq f(\tilde{y}^N) + g(x^N),$$

то следующий результат, позволяет с контролируемой точность восстановить решение задачи (8).

**Теорема 3.** *Пусть нужно решить задачу (8) посредством перехода к задаче (9), исходя из выписанных выше формул. Выбираем в качестве критерия останова УМТ (с $\mu = 0$) выполнение следующих условий*

$$f(\tilde{y}^N) + g(x^N) \leq \varepsilon,\ \|Ax^N - b\|_2 \leq \tilde{\varepsilon}.$$

*Тогда метод гарантированно остановится, сделав не более чем*

$$6 \cdot \max\left\{\sqrt{\frac{LR^2}{\varepsilon}}, \sqrt{\frac{LR}{\tilde{\varepsilon}}}\right\}$$

*итераций, где $R^2 = \|y_*\|_2^2$, $y_*$ – решение задачи (9) (если решение не единственно, то можно считать, что выбирается то $y_*$, которое доставляет минимум $R^2$).*

В основе доказательства теоремы 3 лежит свойство прямо-двойственности УМТ (с $\mu = 0$). С помощью описанной конструкции (и ряда других конструкций) в третьей главе, в частности, изучается задача расчета матрицы корреспонденций и задача поиска равновесия в Смешанной Модели распределения транспортных потоков по ребрам.

В **Главе 4** изучаются задачи оптимизации в пространствах огромных ($\gg 10^7$) размеров (в основном отличные от рассматриваемых ранее "транспортных" задач). Исследуется роль разреженности постановки задачи. Активно используются рандомизированные методы. Изложение построено вокруг классической задачи поиска вектора PageRank. В данной главе, в частности, предложено несколько новых способов решения этой задачи.

Ранее уже отмечалось, что "базисными" методами в диссертации являются МЗС и БГМ. Если БГМ в варианте УМТ уже был описан, то описание МЗС в автореферате пока еще не встречалось (хотя в самой диссертации МЗС параллельно с БГМ активно использовался уже с первой главы). С одной стороны, в плане оценок числа итераций ничего нового получить нельзя – оценки МЗС соответствуют оценкам работы УМТ, когда inf достигается при $\nu = 0$. С другой стороны, МЗС является по структуре намного более простым методом, что позволяет более точно контролировать стоимость итерации, особенно для рандомизированных методов. В данной главе это существенно. Для краткости изложения далее сразу будет описан достаточно общий вариант МЗС для выпуклых задач стохастической условной оптимизации (без ограничения общности $g(x) \leq 0$, см. (10)). В основном в диссертации используется более простой не условный вариант этого метода.

Рассмотрим задачу выпуклой условной оптимизации

$$f(x) \to \min_{g(x) \leq 0,\ x \in Q}. \tag{10}$$

Под решением этой задачи будем понимать такой $\bar{x}^N \in Q \subseteq \mathbb{R}^n$, что с вероятностью $\geq 1 - \sigma$ имеет место неравенство (определение $M_f$, $M_g$ см. ниже)

$$f(\bar{x}^N) - f_* \leq \varepsilon_f = \frac{M_f}{M_g} \varepsilon,\ g(\bar{x}^N) \leq \varepsilon, \tag{11}$$

где $f_* = f(x_*)$ – оптимальное значение функционала в задаче (10), $x_*$ – решение задачи (10). Выберем точку старта $x^1 = \arg\min_{x \in Q} d(x)$. Считаем, что $d(x^1) = 0$, $\nabla d(x^1) = 0$. Как и для УМТ определим $R^2 = V(x_*, x^1)$, где $x_*$ – решение задачи (10) (если решение не единственно, то выбирается то, которое доставляет минимум $V(x_*, x^1)$). Рассмотрим сначала случай ограниченного множества. Тогда можно также определить $\bar{R}^2 = \max_{x, y \in Q} V(y, x)$. Будем считать, что имеется такая последовательность независимых случайных величин $\{\xi^k\}$ и последовательности $\{\nabla_x f(x, \xi^k)\}$, $\{\nabla_x g(x, \xi^k)\}$, $k = 1,...,N$, что имеют место следующие соотношения

$$E_{\xi^k}\left[\nabla_x f(x, \xi^k)\right] = \nabla f(x), \quad E_{\xi^k}\left[\nabla_x g(x, \xi^k)\right] = \nabla g(x); \quad (12)$$

$$\left\|\nabla_x f(x, \xi^k)\right\|_*^2 \leq M_f^2, \quad \left\|\nabla_x g(x, \xi^k)\right\|_*^2 \leq M_g^2 \quad (13)$$

или

$$E_{\xi^k}\left[\left\|\nabla_x f(x, \xi^k)\right\|_*^2\right] \leq M_f^2, \quad E_{\xi^k}\left[\left\|\nabla_x g(x, \xi^k)\right\|_*^2\right] \leq M_g^2. \quad (14)$$

На каждой итерации $k = 1,...,N$ нам доступен стохастический (суб-)градиент $\nabla_x f(x, \xi^k)$ или $\nabla_x g(x, \xi^k)$ в одной, выбранной нами (методом), точке $x^k$.

Определим оператор "проектирования" согласно введенному прокс-расстоянию

$$\text{Mirr}_{x^k}(v) = \arg\min_{y \in Q}\left\{\langle v, y - x^k \rangle + V(y, x^k)\right\}.$$

МЗС для задачи (10) будет иметь вид

### Метод зеркального спуска (МЗС)

$$\boxed{\begin{aligned} x^{k+1} &= \text{Mirr}_{x^k}\left(h_f \nabla_x f(x^k, \xi^k)\right), \quad \text{если } g(x^k) \leq \varepsilon, \\ x^{k+1} &= \text{Mirr}_{x^k}\left(h_g \nabla_x g(x^k, \xi^k)\right), \quad \text{если } g(x^k) > \varepsilon, \end{aligned}} \quad (15)$$

где $h_g = \varepsilon/M_g^2$, $h_f = \varepsilon/(M_f M_g)$, $k = 1,...,N$. Обозначим через $I$ множество индексов $k$, для которых $g(x^k) \le \varepsilon$. Введём также обозначения

$$[N] = \{1,...,N\},\ J = [N]\setminus I,\ N_I = |I|,\ N_J = |J|,\ \bar{x}^N = \frac{1}{N_I}\sum_{k\in I} x^k.$$

В сформулированной далее теореме, предполагается, что последовательность $\{x^k\}_{k=1}^{N+1}$ генерируется методом (15).

**Теорема 4**. *Пусть справедливы условия (12), (14). Тогда при любом*

$$N \ge \frac{2M_g^2 R^2}{\varepsilon^2} + 1$$

*выполняются неравенства $N_I \ge 1$ (с вероятностью $\ge 1/2$) и*

$$E\left[f(\bar{x}^N)\right] - f_* \le \varepsilon_f,\ g(\bar{x}^N) \le \varepsilon.$$

*Пусть справедливы условия (12), (13). Тогда при любом*

$$N \ge \frac{81 M_g^2 \bar{R}^2}{\varepsilon^2}\ln\left(\frac{1}{\sigma}\right)$$

*с вероятностью $\ge 1-\sigma$ выполняются неравенства $N_I \ge 1$ и*

$$f(\bar{x}^N) - f_* \le \varepsilon_f,\ g(\bar{x}^N) \le \varepsilon,$$

*т.е. выполняются неравенства (11).*

В случае неограниченного множества $Q$ при отсутствии ограничения $g(x) \le 0$ приведённая теорема останется верной, если заменить в её формулировке $\bar{R}$ на $\tilde{R} = \sup_{x\in\tilde{Q}}\|x - x_*\|$, где $\tilde{Q} = \{x \in Q: \|x - x_*\|^2 \le 65 R^2 \ln(4N/\sigma)\}$. При этом с вероятностью $\ge 1-\sigma/2$ имеет место включение $\{x^k\}_{k=1}^{N+1} \in \tilde{Q}$.

Описанный МЗС (в рандомизированном, но не условном варианте, с $\|\ \| = \|\ \|_1$ и энтропийной прокс-функцией) лежит в основе ряда подходов поиска вектора PageRank (Назина–Поляка, Юдицкого–Лана–Немировского–Шапиро, Григориадиса–Хачияна, см. таблицу 1 ниже). Примечательно, что для метода Григориадиса–Хачияна рандомизация осуществляется при "проектировании" на симплекс, а не при замене градиента его несмещенной случайной оценкой как в других двух подходах.

Известно, что поисковая система Google была создана в качестве учебного проекта студентов Л. Пейжда и С. Брина из Стэнфордского университета. В 1998 г. Лари Пейждом и Сергеем Брином был предложен специальный способ ранжирования web-страниц. Этот способ, также как и довольно большой класс задач ранжирования, возникающих, например, при вычислении индексов цитирования ученых или журналов, сводится к нахождению левого собственного вектора $p_*$ (нормированного на единицу: $\sum_{k=1}^{n} p_k = 1$, $p_k \geq 0$), отвечающего собственному значению 1 некоторой стохастической (по строкам) матрицы $P = \left\| p_{ij} \right\|_{i,j=1}^{n,n}$, т.е. $p_*$ – решение в классе распределений вероятностей системы $p^T = p^T P$, $n \gg 1$. Предполагаем, что имеется всего один класс сообщающихся состояний, поэтому решение единственно.

Приведем краткое резюме оценок сложности основных алгоритмов поиска вектора PageRank. "Сложность" понимается как количество арифметических операций типа умножения двух чисел, которые достаточно осуществить, чтобы с вероятностью не меньше $1-\sigma$ достичь точности решения $\varepsilon$ по целевому функционалу. Во всех приводимых оценках сложности мы опускаем аддитивное слагаемое $\mathrm{O}(n)$, отвечающее за препроцессинг. Оценки сложности методов из последних пяти строчек приводимой далее Таблицы 1 были получены в диссертационной работе.

**Таблица 1. Сравнение свойств методов решения задачи поиска вектора PageRank**

| Метод | Условие | Сложность | Цель (min) |
|---|---|---|---|
| Назина–Поляка | нет | $\mathrm{O}\left(\dfrac{n \ln(n/\sigma)}{\varepsilon^2}\right)$ | $\left\| P^T p - p \right\|_2^2$ |
| методы Ю.Е. Нестерова | $\bar{S}$ | $\mathrm{O}\left(\dfrac{s^2 \ln n}{\varepsilon^2}\right)$ | $\left\| P^T p - p \right\|_2$ |

| Метод | Условие | Сложность | Норма |
|---|---|---|---|
| вариация алгоритма Юдицкого–Лана–Немировского–Шапиро | нет | $\mathrm{O}\left(\dfrac{n\ln(n/\sigma)}{\varepsilon^2}\right)$ | $\|P^T p - p\|_\infty$ |
| Нестерова–Немировского | G, S | $\dfrac{sn}{\alpha}\ln\left(\dfrac{2}{\varepsilon}\right)$ | $\|p - p_*\|_1$ |
| Поляка–Трембы | S | $\dfrac{2sn}{\varepsilon}$ | $\|P^T p - p\|_1$ |
| Д. Спилмана | G, S | $\mathrm{O}\left(\dfrac{s^2}{\alpha\varepsilon}\ln\left(\dfrac{1}{\varepsilon}\right)\right)$ | $\|p - p_*\|_\infty$ |
| вариация метода Григориадиса–Хачияна из параграфов 1, 2 | $\overline{S}$ | $\mathrm{O}\left(\dfrac{s\ln n\ln(n/\sigma)}{\varepsilon^2}\right)$ | $\|P^T p - p\|_\infty$ |
| MCMC в варианте параграфа 1 (см. также Приложение ниже) | SG | $\mathrm{O}\left(\dfrac{\ln n\ln(n/\sigma)}{\alpha\varepsilon^2}\right)$ | $\|p - p_*\|_2$ |
| вариация метода условного градиента из параграфа 2 | $\overline{S}$ | $\mathrm{O}\left(\dfrac{s^2\ln(n/s^2)}{\varepsilon^2}\right)$ | $\|P^T p - p\|_2$ |
| прямой ACRCD* из главы 5 с $\beta = \dfrac{1}{2}$ | S | $\mathrm{O}\left(\dfrac{sn}{\sqrt{\varepsilon}}\right)$ | $\|Ap - b\|_2^2$ |
| двойственный ACRCD* из главы 5 с $\beta = 1/2$ (см. также главу 3) | S | $\mathrm{O}\left(sn\sqrt{\dfrac{\overline{L}R}{\varepsilon}}\right)$ | $\|Ap - b\|_2$ |

Поясним основные сокращения, использованные в Таблице 1:

- *G-условие* – наличие такой web-страницы (например, страницы, отвечающей самой поисковой системе **G**oogle), на которую можно перейти из любой другой web-страницы с вероятностью не меньшей, чем $\alpha \gg n^{-1}$;
- *S-условие* – из каждой web-страницы в среднем выходит не более $s \ll n$ ссылок на другие, то есть имеет место разреженность матрицы $P$ (**S**parsity); если из каждой web-страницы одновременно выходит и входит не более $s \ll n$ ссылок, то будем говорить об $\overline{S}$-условии;

- *SG-условие* – спектральная щель $\alpha$ (**S**pectral **G**ap) матрицы $P$ удовлетворяет условию $\alpha \gg n^{-1}$, где $\alpha$ – расстояние между максимальным по величине модуля собственным значением (числом Фробениуса–Перрона) матрицы $P$ (равным 1) и модулем следующего (по величине модуля) собственного значения. Если выполняется G-условие, то выполняется и SG-условие с $\alpha$ не меньшим, чем в G-условии.

В таблице 1 также использовались следующие обозначения

$$A = \begin{pmatrix} \left(P^T - I\right) \\ 1 \ \ldots \ 1 \end{pmatrix}, \ b = \begin{pmatrix} 0 \\ \vdots \\ 0 \\ 1 \end{pmatrix}, \ \bar{L} = \left(\frac{1}{n+1}\sum_{k=1}^{n+1}\|A_k\|_2\right)^2,$$

где $A_k$ – вектор, находящийся в $k$-й строке матрицы $A$, а $R$ определяется в теореме 3 для пары задач (8), (9), в которых $x = p$, $g(x) = \frac{1}{2}\|x\|^2$, $Q = \mathbb{R}^n$.

Таблица 1, на наш взгляд, хорошо демонстрирует общую ситуацию: погружаясь в специфику задачи, налагая дополнительные условия, как правило, всегда можно ускориться. В задачах с разреженной структурой наряду с оценками числа итераций метода большую роль начинает играть то, как именно организованы вычисления на каждой итерации. Работу одного и того же метода можно по-разному организовать. В четвертой главе данный тезис демонстрируется на примере перезаписи метода условного градиента (см. таблицу 1). Также таблица 1 и глава 4 в целом демонстрируют и то, что оптимальный метод – совсем не обязательно тот метод, который делает наименьшее число итераций (даже если ограничиться только методами первого порядка). Нужно оценивать общие трудозатраты. К сожалению, в таких понятиях нет никакой теории нижних оценок. Для числа итераций, а точнее числа обращений к локальному оракулу за информацией о функции, такая теория была построена почти 40 лет назад А.С. Немировским. Наконец, из анализа таблицы 1 следует важный практический вывод, что одну и ту же задачу можно по-разному поставить (записать). Другими словами, критерий качества решения можно выбирать различными способами. Иногда удачно выбранный критерий может сильно ускорить поиск решения.

В **Главе 5** анализируются методы порядка ниже первого: покомпонентные методы, спуски по направлению и методы нулевого порядка (безградиентные методы). Нетривиальность заключается в том, что рассматри-ваются рандомизированные и при этом ускоренные методы. Число необходимых итераций (как функция от желаемой точности) для ускоренных покомпонентных методов увеличивается в число раз $n$, равное размерности пространства, по сравнению с классическим БГМ, что и не удивительно, поскольку вместо всех $n \gg 1$ компонент градиента на каждой итерации используется только одна. Также нетривиальность в том, что если полный расчет градиента, скажем, требовал полного умножения разреженной матрицы на вектор – $sn$ операций, то пересчет (важно, что именно пересчет, а не расчет) компоненты градиента и выполнение шага итерации в определенных ситуациях требует всего $s$ операций. Таким образом, увеличение числа итераций в $n$ раз компенсируется уменьшением стоимости одной итерации в $n$ раз (в не разреженном случае оговорка об "определенных ситуациях" существенно ослабляется). Но выгода от использования покомпонентных методов, как правило, есть из-за того, что в таких методах вместо константы Липшица градиента по худшему направлению (максимального собственного значения матрицы Гессе функционала) в оценки числа итераций входит "средняя" константа Липшица, оценивающаяся сверху средним арифметическим суммы диагональных элементов (следа) матрицы Гессе, т.е. средним арифметическим всех собственных чисел матрицы Гессе. Разница в этих константах для матриц Гессе, состоящих из элементов одного порядка, может равняться по порядку $n$. На данный момент известно довольного много примеров применения покомпонентных методов для решения задач огромных размеров, в частности, приложений для задач моделирования сетей больших размеров и анализе данных. В пятой главе покомпонентные методы изучаются сквозь призму именно таких приложений. При этом в данной главе получены новые результаты, позволяющие перенести на покомпонентные методы многие свойства полноградиентных методов: прямо-двойственность, ограниченность итерационной последовательности и т.п.

Ранее в таблице 1 нам уже встречались покомпонентные методы. В частности, если в последней строчке таблицы использовать не покомпонентный метод, а БГМ, то вместо $\bar{L}$ нужно было бы писать $L = \sigma_{\max}(A)$. При этом обычно $L \gg \bar{L}$.

В развитие примера 3 из главы 2, в котором использовалась Minimal Mutual Information Model для восстановления матрицы корреспонденций по замерам потоков на линках (ребрах) в большой компьютерной сети, изучена и другая модель восстановления матрицы корреспонденций (Tomogravity Model), приводящая к задаче

$$\frac{L}{2}\|Ax - b\|_2^2 + \frac{\mu}{2}\|x - x_g\|_2^2 \to \min_{x \in \mathbb{R}^n}.$$

Для этой задачи можно построить двойственную

$$\frac{1}{2\mu}\left(\|x_g - A^T y\|_2^2 - \|x_g\|_2^2\right) + \frac{1}{2L}\left(\|y + b\|_2^2 - \|b\|_2^2\right) \to \min_{y \in \mathbb{R}^m}.$$

В реальных приложениях битовая матрица Кирхгофа $A$ ($m \times n \sim 10^5 \times 10^8$) является сильно разреженной (в среднем $s \ll m$ ненулевых элементов в каждом столбце и в среднем $\tilde{s} = sn/m$ ненулевых элементов в каждой строке). Используя эту специфику, в главе 5 предлагаются ускоренные покомпонентные методы решения прямой и двойственной задачи,[8] общее время работы которых (с точностью до множителя $\sim \ln(\varepsilon^{-1})$ – это отражено волной в $\tilde{O}(\ )$), соответственно, $T^{прям} = \tilde{O}\left(sn\sqrt{Ls/\mu}\right)$, $T^{двойств} = \tilde{O}\left(sn\sqrt{L\tilde{s}/\mu}\right)$. Отсюда можно сделать довольно неожиданный вывод: при $m \ll n$ стоит решать прямую задачу, а в случае $m \gg n$ – двойственную. Первый случай соответствует приложениям к изучению больших (компьютерных) сетей. Второй случай соответствует задачам из области анализа данных.

К сожалению, в общем случае для рандомизаций,[9] отличных от покомпонентных, отмеченные выше преимущества исчезают (во всяком случае, с точки зрения имеющейся сейчас теории). Зато, используя рандомизацию на евклидовой сфере, можно существенно сократить число итераций за счет евклидовой асимметрии множества, на котором происходит оптимизация (эти результаты представляются не только новыми, но и весьма неожиданными). Пусть рассматривается задача минимизации $\mu_p$-сильно ($\mu_p \geq 0$) выпуклой в $p$-норме функции $f(x)$ на множестве $Q$, например, единичном шаре в $p$-норме с оракулом, выдающим производную по случайно выбран-ному на евклидовой сфере направлению или ее дискретную аппроксимацию – для безградиентного оракула. Пусть выбрана норма $\|\ \| = \|\ \|_p$, $1 \le p \le 2$, $1/p + 1/q = 1$, $R^2$ – "расстояние" Брэгмана от точки старта до решения, согласованное с этой нормой (см. главу 2 выше). Пусть в предположении 1 $\|\ \| = \|\ \|_2$ и $M_2$ соответствует $L_0$ ($\nu = 0$), а $L_2$ соответствует $L_1$ ($\nu = 1$). Тогда существуют такие модификации МЗС и БГМ, которые приводят к следующим (оптимальным) оценкам (следует сравнить эти оценки с оценкой (4) в двух наиболее интересных ситуациях $\nu = 0$ и $\nu = 1$), соответственно,

$$\min\left\{\mathrm{O}\!\left(\frac{M_2^2 R^2}{\varepsilon^2}\right)\tilde{\mathrm{O}}\!\left(n^{2/q}\right), \tilde{\mathrm{O}}\!\left(\frac{M_2^2}{\mu_p \varepsilon}\right)\tilde{\mathrm{O}}\!\left(n^{2/q}\right)\right\}, \qquad (16)$$

$$\min\left\{\mathrm{O}\!\left(\sqrt{\frac{L_2 R^2}{\varepsilon}}\right)\tilde{\mathrm{O}}\!\left(n^{1/q+1/2}\right), \tilde{\mathrm{O}}\!\left(\sqrt{\frac{L_2}{\mu_p}}\right)\tilde{\mathrm{O}}\!\left(n^{1/q+1/2}\right)\right\}.$$

---

[8] Для сильно выпуклых постановок прямо-двойственность используемых методов не требуется, поскольку есть сходимость по аргументу и можно восстанавливать решение сопряженной задачи по "модельным" формулам, связывающим прямые переменные с двойственными – такие формулы, как правило, получаются при построении двойственной задачи.

[9] Рандомизация задается способом выбора случайного направления, по которому оракул выдает производную по направлению или её оценку – для методов нулевого порядка используется конечная разность типа Кифера–Вольфовица.

Получены оценки на уровень допустимого шума не случайной природы, при котором выписанные оценки сохраняют вид. Оценки можно распространить (подобно (4), (5)) и на задачи стохастической оптимизации.

В **Главе 6** рассматриваются задачи онлайн-оптимизации с оракулом первого порядка и ниже. Разобрано много конкретных примеров задач оптимизации (оценка для оракула первого порядка соответствует оценке (4) при $\nu = 0$ независимо от степени гладкости и стохастичности). Хотелось бы также отметить методическую находку, в основе которой лежат две онлайн-версии МЗС, позволившую единообразно описать многообразие результатов в этой области. Оригинальным в этой главе является доказательство теоремы, устанавливающей справедливость оценки (16) для задач (стохастической) онлайн-оптимизации (в случае, когда минимум достигается на втором аргументе, формулу (16) удалось доказать только при $q = p = 2$). Выписаны условия на уровень допустимого шума не случайной природы, при которых полученные оценки сохраняют вид.

В **Заключении** приводятся основные результаты и выводы диссертации. Намечаются направления дальнейшего развития. В частности, отмечается перспективность развития концепции конструирования эффективных алгоритмов для структурно сложных задач с помощью небольшого числа базисных методов и развиваемых в диссертации операций над ними.

В **Приложении** описывается формализм современной теории макросистем с точки зрения стохастической химической кинетики. Этот формализм активно используется в первой главе и немного в четвертой главе. Изложим этот формализм с помощью одного из наиболее простых и, одновременно, интересных примеров. В диссертации удалось именно таким образом проинтерпретировать необходимость ранжирования web-страниц согласно вектору PageRank. Насколько нам известно, несмотря на кажущуюся простоту и очевидность, приводимый далее подход является оригинальным.

Имеется $N \gg 1$ пользователей, которые случайно (независимо) блуждают в непрерывном времени по ориентированному графу (на $m$ вершинах – в отличие от главы 4 здесь используется обозначение $m$ для числа вершин графа, а не $n$) с эргодической инфинитезимальной матрицей $\Lambda$ (для этого достаточно, чтобы был всего один класс сообщающихся состояний). Назовем вектор $p$ (из единичного симплекса) PageRank, если $\Lambda p = 0$. Обозначим через $n_i(t)$ – число пользователей на $i$-й странице в момент времени $t \geq 0$. Не сложно показать (теорема Гордона–Ньюэлла), что $n(t)$ асимптотически имеет мультиномиальное распределение с вектором параметров PageRank $p$, т.е.

$$\lim_{t \to \infty} P(n(t) = n) = \frac{N!}{n_1! \cdot \ldots \cdot n_m!} p_1^{n_1} \cdot \ldots \cdot p_m^{n_m}.$$

Следовательно (неравенство Хефдинга в гильбертовом пространстве),[10]

$$\lim_{t \to \infty} P\left( \left\| \frac{n(t)}{N} - p \right\|_2 \geq \frac{2\sqrt{2} + 4\sqrt{\ln(\sigma^{-1})}}{\sqrt{N}} \right) \leq \sigma.$$

Этот же результат можно получить, рассмотрев соответствующую систему унарных химических реакций. Переход одного из пользователей из вершины $i$ в вершину $j$ означает превращение одной молекулы вещества $i$ в одну молекулу вещества $j$, $n_i(t)$ – число молекул $i$-го типа в момент времени $n_i(t)$. Каждое ребро графа соответствует определенной реакции (пре-вращению). Интенсивность реакций определяется матрицей $\Lambda$ и числом молекул, вступающих в реакцию (закон действующих масс). Условие $\Lambda p = 0$ – в точности соответствует условию унитарности в стохастической химической кинетике.

---

[10] На базе этой оценки в главе 4 строится оценка (хорошо параллелизуемого) численного метода MCMC поиска вектора PageRank. Здесь только требуется оговорка, что выписанная в таблице 1 оценка сложности этого метода была установлена при дополнительном предположении (впрочем, его можно ослабить), что вероятности на ребрах (линках), выходящих из одной вершины, одинаковы между собой, и это имеет место для каждой вершины web-графа.

**Публикации автора по теме диссертации в журналах из списка ВАК**

24. Гасников А.В., Лагуновская А.А., Морозова Л.Э. О связи имитационной логит динамики в популяционной теории игр и метода зеркального спуска в онлайн оптимизации на примере задачи выбора кратчайшего маршрута // Труды МФТИ. – 2015. – Т. 7. № 4. – С. 104–113.
25. Гасников А.В., Лагуновская А.А., Усманова И.Н., Федоренко Ф.А. Безградиентные прокс-методы с неточным оракулом для негладких задач выпуклой стохастической оптимизации на симплексе // Автоматика и телемеханика. – 2016. – № 10. – С. 57–77.
26. Гасников А.В., Нестеров Ю.Е., Спокойный В.Г. Об эффективности одного метода рандомизации зеркального спуска в задачах онлайн оптимизации // ЖВМ и МФ. – 2015. – Т. 55. № 4. – С. 582–598.

### Монографии и учебные пособия автора

27. Бузун Н.О., Гасников А.В., Гончаров Ф.О. Горбачев О.Г., Гуз С.А., Крымова Е.А., Натан А.А., Черноусова Е.О. Стохастический анализ в задачах. Учебное пособие. Часть 1. Под ред. А.В. Гасникова. М.: МФТИ, 2016. – 212 с. – Режим доступа: http://arxiv.org/pdf/1508.03461v2.pdf
28. Гасников А.В., Кленов С.Л., Нурминский Е.А., Холодов Я.А., Шамрай Н.Б. Введение в математическое моделирование транспортных потоков. Под ред. А.В. Гасникова. 2-е изд. М.: МЦНМО, 2013. – 427 с. – Режим доступа: http://www.mou.mipt.ru/gasnikov1129.pdf

### Работы автора по теме диссертации

29. Аникин А.С., Гасников А.В., Горнов А.Ю., Двуреченский П.Е., Семенов В.В. Параллелизуемые двойственные методы поиска равновесий в смешанных моделях распределения потоков в больших транспортных сетях. Труды 40-й международной школы-конференции "Информационные технологии и системы – 2016". Россия, Санкт-Петербург (Репино), 25–30 сентября 2016. – 8 с. – Режим доступа: https://arxiv.org/ftp/arxiv/papers/1604/1604.08183.pdf
30. Аникин А.С., Гасников А.В., Двуреченский П.Е., Тюрин А.И., Чернов А.В. Двойственные подходы к задачам минимизации сильно выпуклых функционалов простой структуры при аффинных ограничениях // e-print, 2016. – 16 с. – Режим доступа: https://arxiv.org/ftp/arxiv/papers/1602/1602.01686.pdf
31. Ващенко М.П., Гасников А.В., Молчанов Е.Г., Поспелова Л.Я., Шананин А.А. Вычислимые модели и численные методы для анализа тарифной политики железнодорожных грузоперевозок. М.: ВЦ РАН, 2014. – 52 с. – Режим доступа: https://arxiv.org/ftp/arxiv/papers/1501/1501.02205.pdf
32. Гасников А.В. Марковские модели макросистем // e-print, 2014. – 34 с. – Режим доступа: https://arxiv.org/ftp/arxiv/papers/1412/1412.2720.pdf
33. Гасников А.В., Гасникова Е.В., Мациевский С.В. Прямо-двойственный метод зеркального спуска для условных задач стохастической композит-

Гасников Александр Владимирович

# Эффективные численные методы поиска равновесий в больших транспортных сетях

АВТОРЕФЕРАТ